\documentclass[12pt,reqno]{amsart}
\usepackage{amsmath,amsfonts,amssymb,amsthm,amsxtra,amscd}
\begin{document}

\title[The Horrocks correspondence]{The Horrocks correspondence for 
       coherent sheaves on projective spaces}

\author[I. Coand\u{a}]{Iustin Coand\u{a}}
\address{Institute of Mathematics of the Romanian Academy, P. O. Box 
         1-764, RO-014700, Bucharest, Romania}
\email{Iustin.Coanda@imar.ro}



\thanks{Partially supported by 
CNCSIS grant ID-PCE no.51/28.09.2007 (code 304).}

\begin{abstract}
We establish an equivalence between the stable category of coherent 
sheaves (satisfying a mild restriction) on a projective space and the 
homotopy category of a certain class of minimal complexes of free modules 
over the exterior algebra Koszul dual to the homogeneous coordinate 
algebra of the projective space. We also relate these complexes to the 
Tate resolutions of the respective sheaves. In this way, we extend 
from vector bundles to coherent sheaves the results of G. Trautmann and 
the author \cite{ctr}, which interpret in terms of the BGG correspondence  
the results of Trautmann \cite{trm} about the correspondence of 
Horrocks \cite{ho1}, \cite{ho2}. We also give direct proofs of the 
BGG correspondences for graded modules and for coherent sheaves and 
of the theorem of Eisenbud, Fl\o ystad and Schreyer \cite{efs}
describing the 
linear part of the Tate resolution associated to a coherent sheaf. 
Moreover, we provide an explicit description of the quotient of the 
Tate resolution by its linear strand corresponding to the module of 
global sections of the various twists of the sheaf.  
\end{abstract}

\maketitle

\section*{Introduction}

Two locally free sheaves $E$ and $E^{\prime}$ on the projective space 
${\mathbb P}^n$ 
over a field $k$ are {\it stably equivalent} if there exist finite direct 
sums of invertible sheaves ${\mathcal O}_{\mathbb P}(a)$, $a\in {\mathbb Z}$, 
$L$ and $L^{\prime}$ such that 
$E\oplus L \simeq E^{\prime}\oplus L^{\prime}$. 
Let $S = k[X_0,\dots ,X_n]$ be the homogeneous coordinate ring of 
${\mathbb P}^n$. ${\mathbb P}^n$ being a quotient of $V\setminus \{0\}$, where 
$V = k^{n+1}$, $S$ can be identified with the symmetric algebra $S(V^{\ast})$ 
of the dual vector space $V^{\ast}$. Let $\Lambda := \bigwedge (V)$ be the 
exterior algebra of $V$. For $0 < i < n$, the graded $S$-module 
$\text{H}^i_{\ast}E := \bigoplus_{d\in {\mathbb Z}}\text{H}^i(E(d))$ is an 
invariant for stable equivalence. However, these cohomology $S$-modules alone 
do not determine uniquely the stable equivalence class of $E$. G. Horrocks 
\cite{ho1} showed that the stable equivalence class is determined by these 
modules an by a sequence of extension classes. Unfortunately, the arguments 
of the group $\text{Ext}^1$ in which anyone of these extension classes lives 
depend on the previous extension classes. This inconvenience was removed by 
G. Trautmann \cite{trm} who showed that the stable equivalence class is 
determined by a system of matrices whose entries are (essentially) elements 
of the exterior algebra $\Lambda$. Trautmann's approach is related to the 
approach from Horrocks' paper \cite{ho2}. 

The meaning of the matrices considered by Trautmann was clarified, following a 
suggestion of W. Decker, by Trautmann and the author in \cite{ctr} using the 
Bernstein-Gel'fand-Gel'fand functors. These functors originate in the 
following easy observation: giving a linear complex of graded free 
$S$-modules:  
\[
\cdots \longrightarrow S(p)\otimes_k N_p\longrightarrow 
S(p+1)\otimes_k N_{p+1} \longrightarrow \cdots 
\]
is equivalent to giving a (left) $\Lambda$-module structure on the graded 
$k$-vector space $N := \bigoplus_{p\in {\mathbb Z}}N_p$. One denotes the 
above complex by $\text{F}(N)$. Similarly, to a graded $S$-module $M$ one 
can associate a linear complex $\text{G}(M)$ of graded free $\Lambda$-modules: 
\[
\cdots \longrightarrow M_p\otimes_k{\textstyle \bigwedge}(V^{\ast})(p) 
\longrightarrow M_{p+1}\otimes_k{\textstyle \bigwedge}(V^{\ast})(p+1) 
\longrightarrow \cdots 
\]
where one considers on the exterior algebra $\bigwedge (V^{\ast})$, graded 
such that $V^{\ast}$ has degree $-1$, the structure of left $\Lambda$-module 
defined by contraction. The technical reason for which one uses $\bigwedge 
(V^{\ast})(p)$ instead of $\Lambda (p)$ is that $\text{F}(\bigwedge 
(V^{\ast}))$ is the Koszul resolution of $S/(X_0,\ldots ,X_n)$ 
\[
0\rightarrow S(-n-1)\otimes_k\overset{n+1}{\wedge}V^{\ast}\rightarrow 
\cdots \rightarrow S(-1)\otimes_kV^{\ast}\rightarrow S\rightarrow 0.
\]
The idea of I.N. Bernstein, I.M. Gel'fand and S.I. Gel'fand \cite{bgg} 
was to extend 
these functors to complexes of modules by the formula $\text{F}(N^{\bullet}) 
:= \text{tot}(X^{\bullet \bullet})$, where $N^{\bullet}$ is a complex of 
graded $\Lambda$-modules and $X^{\bullet \bullet}$ is the double complex 
defined by $X^{p,\bullet} := \text{F}(N^p)$, and similarly for $\text{G}
(M^{\bullet})$. 

Now, consider the linear complex $\bigoplus_{i=1}^{n-1}
\text{T}^{-i}\text{G}(\text{H}^i_{\ast}E)$ (T the translation functor for 
complexes) with terms $G^p = \bigoplus_{i=1}^{n-1}\text{H}^i(E(p-i))\otimes_k 
\bigwedge (V^{\ast})(p-i)$ and let $\lambda$ be its differential. 
The first main result of the paper \cite{ctr} asserts that the stable 
equivalence class of $E$ is determined by a perturbation 
$d = \lambda + \delta$ of $\lambda$ obtained by addition of terms of degree 
$\geq 2$. Here ``perturbation'' means that $d \circ d = 0$, i.e., 
$G^{\bullet} := ((G^p)_{p\in {\mathbb Z}},d)$ is a complex, and 
``obtained by addition of terms of degree $\geq 2$'' means that:
\[
{\delta}^p(\text{H}^i(E(p-i))\otimes_k{\textstyle \bigwedge}(V^{\ast})(p-i)) 
\subseteq {\textstyle \bigoplus}_{j<i}\text{H}^j(E(p+1-j))\otimes_k 
{\textstyle \bigwedge}(V^{\ast})(p+1-j).
\]
The second main result of \cite{ctr} relates the Horrocks correspondence 
to the BGG correspondence via the results of Eisenbud, Fl\o ystad and 
Schreyer \cite{efs} about Tate resolutions over the exterior algebra. Let 
$\mathcal F$ be a coherent sheaf on ${\mathbb P}^n$. Eisenbud et al. 
\cite{efs} show that there is a unique (up to isomorphism) perturbation 
of the differential of the linear complex $\bigoplus_{i=0}^n\text{T}^{-i} 
\text{G}(\text{H}^i_{\ast}{\mathcal F})$ 
obtained by addition of terms of degree 
$\geq 2$ such that the resulting complex $I^{\bullet}$ is {\it acyclic}. 
It is shown in \cite{ctr} that the complex $G^{\bullet}$ 
which determines the stable equivalence class of a locally free sheaf $E$ 
can be obtained from its {\it Tate resolution} $I^{\bullet}$ by removing 
the linear strands $\text{G}(\text{H}^0_{\ast}E)$ and $\text{T}^{-n}
\text{G}(\text{H}^n_{\ast}E)$. 

In this paper, we generalize the results from \cite{ctr} to the case of 
coherent sheaves using different, more natural, arguments: while in 
\cite{ctr} one avoids the use of the BGG correspondence, the proofs in the 
present paper depend on it. In the first section we show that, using 
arguments close to the arguments of Horrocks \cite{ho1}, one can extend 
from vector bundles to coherent sheaves the splitting criterion of 
Horrocks and his criterion of stable equivalence. The extension to coherent 
sheaves of the first criterion is a result  
obtained recently by Abe and Yoshinaga 
\cite{abe} (see, also, Bertone and Roggero \cite{ber}). 

In the second section we introduce and prove the properties of the 
BGG functors needed in the proof of the Horrocks correspondence. 
These are : (1) the BGG equivalence between the bounded derived category 
of finitely generated graded $S$-modules and the corresponding category 
of $\Lambda$-modules, for which we provide a direct proof, avoiding the 
use of Koszul duality; (2) the easy half of the Koszul duality phenomenon 
which says that if $N$ is a graded $\Lambda$-module then $\text{G}
\text{F}(N)$ is a right resolution  
of $N$ with graded free $\Lambda$-modules; 
(3) using, additionally, the functors $\text{Hom}_S(-,S)$ for $S$-modules 
and $\text{Hom}_k(-,k)$ for $\Lambda$-modules, one deduces, modulo some 
unpleasant sign problems, a functorial (left) free resolution for every 
$\Lambda$-module $N$; (4) a key technical point of the paper of Eisenbud 
et al. \cite{efs} describing the linear part of the  
minimal complex associated to a complex 
of free modules of the form $\text{F}(N^{\bullet})$ or 
$\text{G}(M^{\bullet})$. The last result is a consequence of a general 
lemma about double complexes, see \cite{efs}, (3.5). We explain, in 
Appendix A, that this lemma is a particular case of a general lemma 
well-known in homotopy theory under the name of 
Basic Perturbation Lemma. 

In the third section we establish the Horrocks correspondence for 
coherent sheaves. It asserts that the stable category of coherent sheaves 
$\mathcal F$ on ${\mathbb P}^n$ with the property that $\text{H}^0
{\mathcal F}(-t) = 0$ for $t >> 0$ is equivalent to the homotopy category 
of {\it minimal} complexes $G^{\bullet}$ 
of graded free $\Lambda$-modules whose 
linear part is of the form $\bigoplus_{i=1}^{n-1}\text{T}^{-i}\text{G}(H^i)$, 
where $H^i$ is the $k$-vector space graded dual of a finitely generated 
graded  $S$-module of Krull dimension $\leq i+1$. This is equivalent to the 
fact that $G^{\bullet}$ is minimal and satisfies 
the following three conditions: 

(i) $G^{\bullet}$ \textit{is right bounded and} 
$\text{H}^p(G^{\bullet}) = 0$ \textit{for} $p << 0$,\\
\hspace*{3mm} (ii) $\forall p\in {\mathbb Z}$, $G^p$ 
\textit{is of the form} $\bigoplus_{i=1}^{n-1}\bigwedge (V^{\ast})
(p-i)^{\textstyle c_{pi}}$,\\
\hspace*{3mm} (iii) $\underset{p\rightarrow \infty}{\text{lim}}
(c_{-p,i}/p^{i+1}) = 0$, $i=1,\ldots ,n-1$. 

In the fourth section we relate the Horrocks correspondence and the 
BGG correspondence. We first give a direct proof of the BGG 
description of the bounded derived category of coherent sheaves on 
${\mathbb P}^n$. This proof is based on an elementary comparison lemma 
which is discussed in Appendix B. 
Using the comparison lemma we also get a quick proof of the theorem of 
Eisenbud, Fl\o ystad and Schreyer \cite{efs}, Theorem 4.1., about Tate 
resolutions of coherent sheaves on ${\mathbb P}^n$. Moreover, we provide  
a concrete description of the quotient $I^{\bullet}/\text{G}
(\text{H}^0_{\ast}{\mathcal F})$, where $I^{\bullet}$ is the Tate 
resolution of a coherent sheaf $\mathcal F$ with $\text{H}^0
({\mathcal F}(-t)) = 0$ for $t >> 0$. Using this concrete description 
we derive that the complex $G^{\bullet}$ associated to $\mathcal F$ by 
the Horrocks correspondence can be obtained from its Tate resolution 
$I^{\bullet}$ by removing the linear strands $\text{G}(\text{H}^0_{\ast} 
{\mathcal F})$ and $\text{T}^{-n}\text{G}(\text{H}^n_{\ast}{\mathcal F})$.
    
\vskip3mm
{\bf Notation.}\quad 
Throughout this paper, $V$ will denote  
an $(n+1)$-dimensional vector space over a field $k$, $e_0,\dots ,
e_n$ a fixed basis of $V$ and $X_0,\dots ,X_n$ the dual basis of 
$V^{\ast} := \text{Hom}_k(V,k)$. 

(i) Let $S=S(V^{\ast})=\bigoplus_{i\geq 0}S^i(V^{\ast})
\simeq k[X_0,\dots ,X_n]$ be the symmetric algebra of 
$V^{\ast}$ and $S_+=\bigoplus_{i\geq 1}S^i(V^{\ast})$ 
its irrelevant homogeneous 
ideal. We denote by $S$-Mod the category of graded $S$-modules {\it with all 
the homogeneous components finite dimensional vector spaces}, and with 
morphisms of degree 0. $S$-mod denotes the full subcategory of $S$-Mod 
consisting of finitely generated graded $S$-modules, and $\mathcal P$ denotes 
the full subcategory of $S$-mod consisting of its free objects. 

(ii) If $M$ is an object of $S$-mod, the finitely generated $S$-module 
$M^{\vee}:=\text{Hom}_S(M,S)$ has a natural grading given by 
$(M^{\vee})_d=\text{Hom}_{S\text{-mod}}(M,S(d))$. 

If $M$ is an object of $S$-Mod, the graded dual vector space 
$M^{\ast}:=\bigoplus_{d\in {\mathbb Z}}\text{Hom}_k(M_{-d},k)$ 
has a natural structure of graded $S$-module. 

(iii) Let $\Lambda = \bigwedge (V)=\bigoplus_{i=0}^{n+1}\overset{i}{\wedge}V$ 
be the (positively graded) exterior algebra of $V$. ${\Lambda}_+:=
\bigoplus_{i=1}^{n+1}\overset{i}{\wedge}V$ 
is an ideal of $\Lambda $. Let $\underline{k}$ denote 
the quotient $\Lambda /{\Lambda}_+$. We denote by $\Lambda $-mod the category 
of finitely generated graded left $\Lambda $-modules with morphisms of degree 
0, and by $\mathcal I$ its full subcategory consisting of free objects. 

If $N\in \text{Ob}(\Lambda\text{-mod})$, $\text{soc}(N)$ denotes the 
submodule of $N$ consisting of the elements annihilated by ${\Lambda}_+$. It 
can be identified with $\text{Hom}_{\Lambda}(\underline{k},N)$. 
Remark that $\text{soc}(\Lambda) = \overset{n+1}{\wedge}V$.  

(iv) Let ${\mathbb P} = {\mathbb P}^n = {\mathbb P}(V)$ denote the 
(classical) projective space parametrizing the 1-dimensional vector 
subspaces of $V$. The homogeneous coordinate ring of $\mathbb P$ is $S$. We 
denote by $\text{Coh}\, {\mathbb P}$ the category of coherent sheaves on 
$\mathbb P$ and by $(-)\sptilde : S\text{-mod}\rightarrow \text{Coh}\, 
{\mathbb P}$ the functor of Serre \cite{ser}
associating to a graded $S$-module its sheafification. 
If $\mathcal F$ is a coherent sheaf on $\mathbb P$ and 
$0\leq i\leq n$ we shall denote by $\text{H}^i_{\ast}{\mathcal F}$ the 
graded $S$-module $\bigoplus_{d\in {\mathbb Z}}\text{H}^i({\mathcal F}(d))$.  

(v) We denote by $\text{C}(\mathcal A)$, $\text{C}^b(\mathcal A)$, 
$\text{C}^{\pm}(\mathcal A)$ the categories of complexes in an abelian 
category $\mathcal A$, by $\text{K}(\mathcal A)$, $\text{K}^b(\mathcal A)$, 
$\text{K}^{\pm}(\mathcal A)$ the corresponding homotopy categories, and by 
$\text{D}(\mathcal A)$, $\text{D}^b(\mathcal A)$, $\text{D}^{\pm}(\mathcal A)$ 
the corresponding derived categories. ``T'' will denote the {\it translation 
functor} and ``Con'' the {\it mapping cone}. When $\mathcal A$ is $S$-mod or 
$\Lambda $-mod or $\text{Coh}\, {\mathbb P}$ we shall use the shorter 
notation $\text{C}(S)$, $\text{C}(\Lambda )$, $\text{C}(\mathbb P)$ etc. 

Our main reference for category theory will be Chapter I of the book of 
Kashiwara and Schapira \cite{ks1}. One may also use the books of Kashiwara 
and Schapira \cite{ks2} or Gel'fand and Manin \cite{gma}.

\section{Two criteria of Horrocks}

In this section we include proofs of two introductory results 
which extend to coherent sheaves the 
splitting criterion of Horrocks for vector bundles on projective spaces and 
his criterion characterizing stable equivalences in the same context.  
There are (at least) two recently published proofs of the first result in 
Abe and Yoshinaga \cite{abe} and Bertone and Roggero \cite{ber}. 
We follow, however, Horrocks' original approach. It is based on 
the next theorem, which is usually proved using local cohomology and local 
duality (see, for example, \cite{eis}, (A.4.1.) and (A.4.2)). 
We shall give, for the reader's convenience,  
a direct proof avoiding the use of local cohomology. 

\vskip3mm
{\bf 1.1. Theorem.} (Graded Serre Duality)\quad \textit{If} $M$ \textit{is a 
finitely generated graded} $S$-\textit{module then there exist an exact 
sequence}:
\[
0\rightarrow \text{Ext}^{n+1}_S(M,{\omega}_S)^{\ast}\rightarrow M\rightarrow 
\text{H}^0_{\ast}{\widetilde M}\rightarrow \text{Ext}^n_S(M,{\omega}_S)^{\ast} 
\rightarrow 0
\]
\textit{and isomorphisms}: $\text{H}^i_{\ast}{\widetilde M} \simeq 
\text{Ext}^{n-i}_S(M,{\omega}_S)^{\ast}$, $1 \leq i \leq n$, \textit{where} 
${\omega}_S := S(-n-1)\otimes_k\overset{n+1}{\wedge}V^{\ast}$.
\vskip3mm

\begin{proof}
${\widetilde \omega}_S = {\mathcal O}_{\mathbb P}(-n-1)\otimes_k
\overset{n+1}{\wedge}V^{\ast} \simeq {\omega}_{\mathbb P}$. One knows that 
$\text{H}^n({\omega}_{\mathbb P}) \simeq k$ and that, $\forall a\in 
{\mathbb Z}$: 
\[
\text{Hom}_{{\mathcal O}_{\mathbb P}}({\mathcal O}_{\mathbb P}(-a), 
{\omega}_{\mathbb P})\overset{\sim}{\longrightarrow} 
\text{Hom}_k(\text{H}^n({\mathcal O}_{\mathbb P}(-a)),\text{H}^n
({\omega}_{\mathbb P})).
\]
It follows that if $L$ is a free graded $S$-module of finite rank then there 
exists a functorial isomorphism: 
\[
\text{Hom}_S(L,{\omega}_S)\overset{\sim}{\longrightarrow}(\text{H}^n_{\ast}
{\widetilde L})^{\ast}.
\]
Now, let $0\rightarrow L^{-n-1}\rightarrow \cdots \rightarrow L^0 \rightarrow  
M \rightarrow 0$ be a free resolution of $M$ in $S$-mod. Let $C^{-i}$ be the 
cokernel of $L^{-i-1} \rightarrow L^{-i}$. One has short exact sequences:
\[
0\rightarrow C^{-i-1}\rightarrow L^{-i}\rightarrow C^{-i}\rightarrow 0.
\tag{1}
\]    
We consider the complex $\text{H}^n_{\ast}{\widetilde L}^{\bullet} \simeq 
\text{Hom}_S(L^{\bullet},{\omega}_S)^{\ast}$. Since 
$\text{H}^n_{\ast}$ is right exact, we have exact sequences:
\[
\text{H}^n_{\ast}{\widetilde L}^{-i-1}\rightarrow \text{H}^n_{\ast}
{\widetilde L}^{-i}\rightarrow \text{H}^n_{\ast}{\widetilde C}^{-i}
\rightarrow 0
\]
hence $\text{H}^0(\text{H}^n_{\ast}{\widetilde L}^{\bullet}) \simeq 
\text{H}^n_{\ast}{\widetilde M}$ and $\text{H}^{-i}(\text{H}^n_{\ast}
{\widetilde L}^{\bullet}) \simeq \text{Ker}(\text{H}^n_{\ast}
{\widetilde C}^{-i} \rightarrow \text{H}^n_{\ast}{\widetilde L}^{-i+1})$ 
for $i \geq 1$. Since $\text{H}^p_{\ast}{\widetilde L}^{-j} = 0$ for 
$0 < p < n$, $\forall j$, one deduces easily, using the sheafifications of 
the exact sequences (1), that:
\[
\text{H}^{-i}(\text{H}^n_{\ast}{\widetilde L}^{\bullet})\simeq 
\text{H}^{n-i}_{\ast}{\widetilde M},\  \  0\leq i\leq n-1
\]
and that one has exact sequences:
\begin{gather*}
L^0=\text{H}^0_{\ast}{\widetilde L}^0\rightarrow \text{H}^0_{\ast}
{\widetilde M}\rightarrow \text{H}^{-n}(\text{H}^n_{\ast}
{\widetilde L}^{\bullet})\rightarrow 0\\
L^{-1}=\text{H}^0_{\ast}{\widetilde L}^{-1}\rightarrow \text{H}^0_{\ast}
{\widetilde C}^{-1}\rightarrow \text{H}^{-n-1}(\text{H}^n_{\ast}
{\widetilde L}^{\bullet})\rightarrow 0,
\end{gather*}
hence:
\begin{gather*}
\text{H}^{-n}(\text{H}^n_{\ast}{\widetilde L}^{\bullet})\simeq 
\text{Coker}(M\rightarrow \text{H}^0_{\ast}{\widetilde M})\\
\text{H}^{-n-1}(\text{H}^n_{\ast}{\widetilde L}^{\bullet})\simeq 
\text{Coker}(C^{-1}\rightarrow \text{H}^0_{\ast}{\widetilde C}^{-1}).
\end{gather*}
Finally, applying the Snake Lemma to the diagram:
\[
\begin{CD}
0 @>>> C^{-1} @>>> L^0 @>>> M @>>> 0\\
@. @VVV @VV{\wr}V @VVV\\
0 @>>> \text{H}^0_{\ast}{\widetilde C}^{-1} @>>> \text{H}^0_{\ast}
{\widetilde L}^0 @>>> \text{H}^0_{\ast}{\widetilde M} 
\end{CD}
\]
one gets that: $\text{Coker}(C^{-1}\rightarrow \text{H}^0_{\ast}
{\widetilde C}^{-1}) \simeq \text{Ker}(M\rightarrow \text{H}^0_{\ast}
{\widetilde M})$.  
\end{proof}
\vskip3mm

{\bf 1.2. Lemma.}\quad 
\textit{If} $M\neq 0$ \textit{is a finitely generated graded} 
$S$-\textit{module of projective dimension} $m$ \textit{then} 
$\text{Ext}^m_S(M,S) \neq 0$.
\vskip3mm

\begin{proof}
Let $0\rightarrow L^{-m}\rightarrow \cdots \rightarrow L^0\rightarrow M 
\rightarrow 0$ be a free resolution of $M$ in $S$-mod. If 
$\text{Ext}^m_S(M,S) = 0$ then $L^{-m+1 \vee} \rightarrow L^{-m \vee}$ is 
surjective, hence its kernel $L^{\prime -m+1}$ is free and the sequence: 
\[
0\rightarrow L^{\prime -m+1}\rightarrow L^{-m+1 \vee}\rightarrow 
L^{-m \vee}\rightarrow 0
\]
is split exact. It follows that the dual sequence:
\[
0\rightarrow L^{-m}\rightarrow L^{-m+1}\rightarrow (L^{\prime -m+1})^{\vee} 
\rightarrow 0
\]
is exact. One gets an exact sequence: 
\[
0\rightarrow (L^{\prime -m+1})^{\vee}\rightarrow L^{-m+2}\rightarrow \cdots 
\rightarrow L^0\rightarrow M\rightarrow 0
\]
from which we deduce that the projective dimension of $M$ is $\leq m-1$, a 
contradiction. 
\end{proof}
\vskip3mm

{\bf 1.3. Theorem.} (Horrocks' splitting criterion)\quad 
\textit {Let} $\mathcal F$ \textit{be a coherent sheaf on} ${\mathbb P}^n$ 
\textit{with the property that} $\text{H}^0({\mathcal F}(-t)) = 0$ 
\textit{for} $t >> 0$. \textit{If} $\text{H}^i_{\ast}{\mathcal F} = 0$ 
\textit{for} $0 < i < n$ \textit{then} $\mathcal F$ \textit{is a direct 
sum of invertible sheaves} ${\mathcal O}_{\mathbb P}(a)$, $a\in {\mathbb Z}$. 
\vskip3mm

\begin{proof}
By hypothesis, $M := \text{H}^0_{\ast}{\mathcal F}$ is a {\it finitely 
generated} graded $S$-module. It follows that $M \rightarrow \text{H}^0_{\ast} 
{\widetilde M}$ is an isomorphism. One deduces, now, from the hypothesis and 
from Theorem 1.1., that $\text{Ext}^i_S(M,{\omega}_S) = 0$, $\forall i > 0$. 
It follows, from Lemma 1.2., that $M$ is a graded free $S$-module, hence a 
direct sum of graded $S$-modules of the form $S(a)$, $a\in {\mathbb Z}$.
\end{proof}
\vskip3mm

{\bf 1.4. Theorem.} (Horrocks' criterion of stable equivalence)\quad 
\textit{Let} $\phi : {\mathcal F}\rightarrow {\mathcal G}$ 
\textit{be a morphism of coherent sheaves on} 
${\mathbb P}^n$, $n\geq 2$, \textit{with the property that} 
$\text{H}^0\phi (-t)$ \textit{is an isomorphism for} $t >> 0$. \textit{If} 
$\text{H}^i_{\ast}\phi$ \textit{is an isomorphism for} $0 < i < n$ 
\textit{then} $\phi$ \textit{factorizes as}:
\[
{\mathcal F}\hookrightarrow {\mathcal F}\oplus {\mathcal A}
\overset{\sim}{\longrightarrow} {\mathcal G}\oplus {\mathcal B} 
\twoheadrightarrow {\mathcal G} 
\]
\textit{where the first morphism is the canonical inclusion}, $\mathcal A$ 
\textit{and} $\mathcal B$ \textit{are finite direct sums of invertible 
sheaves} ${\mathcal O}_{\mathbb P}(a)$, $a\in {\mathbb Z}$, \textit{and the 
last morphism is the canonical projection}.
\vskip3mm 

\begin{proof}
Choose $m\in {\mathbb Z}$ such that $\text{H}^0\phi (-t)$ is an isomorphism 
for $t > m$ and let $M := \bigoplus_{j\geq -m}\text{H}^0({\mathcal F}(j))$, 
$N := \bigoplus_{j\geq -m}\text{H}^0({\mathcal G}(j))$. Choose an 
epimorphism $g : A \rightarrow N$, with $A$ a finitely generated graded 
free $S$-module. Let $\pi : {\mathcal F}\oplus {\widetilde A} \rightarrow 
{\mathcal G}$ be the epimorphism defined by $\phi$ and $\widetilde g$ and let 
$\mathcal B$ be the kernel of $\pi$. Using the exact sequence:
\[
0\rightarrow {\mathcal B}\rightarrow {\mathcal F}\oplus {\widetilde A}
\rightarrow {\mathcal G}\rightarrow 0
\tag{*}
\]
one sees that $\mathcal B$ satisfies the hypothesis of Theorem 1.3., hence 
$\mathcal B$ is a direct sum of invertible sheaves 
${\mathcal O}_{\mathbb P}(b)$, $b\in {\mathbb Z}$, and, consequently, 
$B := \text{H}^0_{\ast}{\mathcal B}$ is a graded free $S$-module. Applying 
$\text{H}^0_{\ast}$ to the exact sequence (*) and cancellating the 
isomorphism $\bigoplus_{j<-m}\text{H}^0{\mathcal F}(j) 
\overset{\sim}{\rightarrow} \bigoplus_{j<-m}\text{H}^0{\mathcal G}(j)$ one 
gets a short exact sequence:
\[
0\rightarrow B\rightarrow M\oplus A\rightarrow N\rightarrow 0.
\tag{**}
\]
Since $\text{H}^{n-1}_{\ast}\phi$ is an isomorphism, it follows from 
Theorem 1.1. that: 
\[
\text{Ext}^1_S(N,B) \longrightarrow 
\text{Ext}^1_S(M\oplus A,B)
\] 
is an isomorphism, hence $\text{Hom}_S(M\oplus A, 
B) \rightarrow \text{Hom}_S(B,B)$ is surjective, hence the exact sequence 
(**) splits.
\end{proof} 
\vskip3mm 

\section{The BGG functors}
 
{\bf 2.1. Definition.}\quad 
When dealing with the category $\Lambda $-mod one encounters sign problems. 
In order to avoid any complication we shall observe strictly the {\it Koszul 
sign convention} (when two homogeneous symbols $\xi $ and $\eta $ are 
permuted the result is multiplied by $(-1)^{\text{deg}\, \xi \, \cdot \,  
\text{deg}\, \eta}$). 

(i) If $K$, $N\in \text{Ob}(\Lambda \text{-mod})$ the graded $k$-vector 
space $K\otimes_kN$ has a structure of left $\Lambda $-module given by: 
\[
v\cdot (x\otimes y):=(v\cdot x)\otimes y+(-1)^{\text{deg}\, x}x\otimes 
(v\cdot y),\  \  \text{for}\  v\in V.
\]
In particular, we put, for $a\in {\mathbb Z}$, $N(a):=\underline{k}(a)
\otimes_kN$. The grading of $N(a)$ is given by $N(a)_p=N_{p+a}$ and the 
$\Lambda $-module structure by: $(v\cdot y)_{N(a)}=(-1)^a(v\cdot y)_N$, for 
$v\in V$, $y\in N$. With this definition, if $v\in V$ then the morphism of 
$k$-vector spaces $(v\cdot -)_N:N\rightarrow N$ is a morphism in 
$\Lambda $-mod: $N(a)\rightarrow N(a+1)$, $\forall a\in {\mathbb Z}$. If 
$\phi :K\rightarrow N$ is a morphism in $\Lambda $-mod, $\phi (a):K(a)
\rightarrow N(a)$ is just $\phi $ if one forgets the gradings. However, 
if $N^{\bullet}\in \text{Ob}\, \text{C}(\Lambda \text{-mod})$ then 
$N^{\bullet}(a)$ is, by definition, the complex with terms $(N^p(a))_{p\in 
{\mathbb Z}}$ but with $d_{N(a)}:=(-1)^ad_N$ (the differential of a complex 
is a symbol of degree 1!).

(ii) If $K$, $N\in \text{Ob}(\Lambda \text{-mod})$ the graded 
$k$-vector space $\text{Hom}_k(N,K)$ has a structure of left 
$\Lambda $-module given by:
\[
(v\cdot \phi)(y):=v\cdot \phi (y)-(-1)^{\text{deg}\, \phi}\phi (v\cdot y), 
\  \  \text{for}\  v\in V.
\]
In particular, for $K=\underline{k}$, one puts $N^{\ast}:=\text{Hom}_k(N,
\underline{k})$. One has $(N^{\ast})_p=(N_{-p})^{\ast}$ and, for $v\in V$, 
$(v\cdot -)_{N^{\ast}}:(N^{\ast})_p\rightarrow (N^{\ast})_{p+1}$ is 
$(-1)^{p+1}\cdot \, \text{the dual of}\  (v\cdot -)_N: N_{-p-1}\rightarrow 
N_{-p}$. 

The map $\mu :N\rightarrow N^{\ast \ast}$, $\mu (y)(\phi) :=
(-1)^{\text{deg}\, y\, \cdot \, \text{deg}\, \phi}\phi (y)$ (i.e., with 
${\mu}_p:=(-1)^p\text{can}:N_p\rightarrow (N_p)^{\ast \ast}$, $p\in 
{\mathbb Z}$) is an isomorphism in $\Lambda $-mod.

(iii) The map $\alpha :K^{\ast}\otimes_kN^{\ast}\rightarrow 
(K\otimes_kN)^{\ast}$ given by:
\[
\alpha (f\otimes g)(x\otimes y):=(-1)^{\text{deg}\, g\, \cdot \,  
\text{deg}\, x}f(x)g(y)
\]
is an isomorphism in $\Lambda $-mod. In particular, for $a\in {\mathbb Z}$,  
one gets an isomorphism in $\Lambda $-mod $\alpha :N^{\ast}(-a)
\overset{\sim}{\rightarrow} (N(a))^{\ast}$ with 
${\alpha}_p=(-1)^{(p-a)a}\, \text{id}_{(N_{a-p})^{\ast}}$, $p\in {\mathbb Z}$. 

Under these identifications, if $v\in V$ and $a\in {\mathbb Z}$, 
the dual of the morphism $(v\cdot -)_N:N(-a-1)\rightarrow N(-a)$ is 
identified with the morphism $(-1)^a(v\cdot -)_{N^{\ast}}:N^{\ast}(a)
\rightarrow N^{\ast}(a+1)$. 

(iv) We endow 
$\bigwedge (V^{\ast})=\bigoplus_{i=0}^{n+1}\overset{i}{\wedge} V^{\ast}$,  
graded such that $\overset{i}{\wedge} V^{\ast}=\bigwedge (V^{\ast})_{-i}$,  
with the structure of graded left $\Lambda $-module given by: 
\[
v\cdot f_1\wedge \ldots \wedge f_p := \sum_{i=1}^p(-1)^{i-1}
f_i(v)f_1\wedge \ldots \wedge 
\widehat{f_i} \wedge \ldots \wedge f_p,\  \  v\in V,\  f_1,\ldots ,f_p\in 
V^{\ast}.
\]
The unique morphism of left $\Lambda $-modules $\Lambda \rightarrow 
\bigwedge (V^{\ast})^{\ast}$ (resp., $\Lambda \rightarrow \bigwedge (V^{\ast}) 
(-n-1)\otimes_k\overset{n+1}{\wedge} V$) mapping $1\in {\Lambda}_0$ to 
$1\in (\bigwedge (V^{\ast})^{\ast})_0$ 
(resp., to $\text{id}_{\overset{n+1}{\wedge} V} 
\in (\bigwedge (V^{\ast})(-n-1)\otimes_k\overset{n+1}{\wedge} V)_0$) is an 
isomorphism in $\Lambda $-mod. 
\vskip3mm 
  
The following lemma, whose standard proof can be found, for example, in 
\cite{coa}, (4)(i), shows, in particular, that, $\forall a\in {\mathbb Z}$, 
$\bigwedge (V^{\ast})(a)$ is an injective object of $\Lambda$-mod. 

\vskip3mm
{\bf 2.2. Lemma.}\quad 
$\forall N\in \text{Ob}(\Lambda \text{-mod})$, $\forall a\in {\mathbb Z}$, 
\textit{the map}:
\[
\text{Hom}_{\Lambda \text{-mod}}(N,{\textstyle \bigwedge} (V^{\ast})(a))
\longrightarrow 
(N_{-a})^{\ast},\  \  \phi \mapsto {\phi}_{-a} 
\]
\textit{is bijective}.
\vskip3mm

{\bf 2.3. Remark.}\quad
${\phi}_{-a-1}: N_{-a-1}\rightarrow \bigwedge 
(V^{\ast})(a)_{-a-1}=V^{\ast}$ \textit{can be described by}:
\[
{\phi}_{-a-1}(y)(v)=(-1)^a{\phi}_{-a}(v\cdot y),\  \forall y\in N_{-a-1},\  
\forall v\in V,
\]
\textit{or, equivalently, by}: ${\phi}_{-a-1}=(-1)^a\sum_{i=0}^n
({\phi}_{-a}\circ (e_i\cdot -)_N)\otimes X_i$.
\vskip3mm

\begin{proof}
Recalling the Definition 2.1.(i), (ii), one has, for $v\in V$, $\lambda \in 
V^{\ast}$:
\[
(v\cdot \lambda)_{\bigwedge (V^{\ast})(a)}=(-1)^a(v\cdot \lambda)_
{\bigwedge (V^{\ast})}=(-1)^a\lambda (v).
\]
In particular, for $\lambda ={\phi}_{-a-1}(y)$,  
since $\phi$ is a morphism in $\Lambda$-mod:
\[
(-1)^a{\phi}_{-a-1}(y)(v)=(v\cdot {\phi}_{-a-1}(y))_{\bigwedge (V^{\ast})(a)} 
={\phi}_{-a}(v\cdot y).
\]
\end{proof}

{\bf 2.4. Definition.} (The BGG functors)\\
\hspace*{3mm} (i) One defines a functor $\text{F}: \Lambda \text{-mod} 
\rightarrow \text{C}^b(\mathcal P)$ by: 
$
\text{F}(N)^p:=S(p)\otimes_kN_p,\  \  d_{\text{F}(N)}:=\sum_{i=0}^n(X_i\cdot 
-)_S\otimes (e_i\cdot -)_N.
$
F can be extended to a functor $\text{F}: \text{C}^b(\Lambda \text{-mod})
\rightarrow \text{C}^b(\mathcal P)$ by putting $\text{F}(N^{\bullet}):=
\text{tot}(X^{\bullet \bullet})$, where $X^{\bullet \bullet}$ is the double 
complex with $X^{p,\bullet}:=\text{F}(N^p)$ and with $d_X^{\prime p}:
X^{p, \bullet}\rightarrow X^{p+1, \bullet}$ equal to $\text{F}(d^p_N)$. 

(ii) One defines a functor $\text{G}: S\text{-Mod}\rightarrow \text{C}
(\mathcal I)$ by: 
$
\text{G}(M)^p:=M_p\otimes_k{\textstyle \bigwedge}(V^{\ast})(p),\  \  
d_{\text{G}(M)}:= \sum_{i=0}^n(X_i\cdot -)_M\otimes (e_i\cdot -)_
{\bigwedge (V^{\ast})}.
$
G can be extended, in a similar way,  
to a functor $\text{G}: \text{C}^b(S\text{-Mod})
\rightarrow \text{C}(\mathcal I)$. The (extended) functor G maps 
$\text{C}^b(S\text{-mod})$ to $\text{C}^+(\mathcal I)$. 

(iii) F {\it and} G {\it commute with the translation functors and with 
mapping cones}. 

\vskip3mm
{\bf 2.5. Definition.}\quad 
Let $\Phi : {\mathcal A}^\text{op}\rightarrow {\mathcal B}$ be an additive 
contravariant functor between two additive categories $\mathcal A$ and 
$\mathcal B$. If $X^{\bullet}\in \text{Ob}\, \text{C}(\mathcal A)$, 
one defines 
a complex $\Phi (X^{\bullet})\in \text{Ob}\, \text{C}(\mathcal B)$ by: 
\[
\Phi (X^{\bullet})^p:=\Phi (X^{-p}),\  \  d^p_{\Phi (X)}:=(-1)^{p+1}\Phi 
(d_X^{-p-1}): \Phi (X^{-p})\rightarrow \Phi (X^{-p-1}).
\]
\hspace*{3mm} {\it For example}, 
if $M^{\bullet}\in \text{Ob}\, \text{C}(S\text{-mod})$ one can 
define the complex $M^{\bullet \vee}\in \text{Ob}\, \text{C}(S\text{-mod})$ 
and 
if $M^{\bullet}\in \text{Ob}\, \text{C}(S\text{-Mod})$ (resp., $N^{\bullet}\in 
\text{Ob}\, \text{C}(\Lambda \text{-mod})$) one can define the complex 
$M^{\bullet \ast}\in \text{Ob}\, \text{C}(S\text{-Mod})$ (resp., 
$N^{\bullet \ast} \in \text{Ob}\, \text{C}(\Lambda \text{-mod})$).

Furthermore, if $X^{\bullet \bullet}$ is a double complex in $\mathcal A$ one 
defines a double complex $\Phi (X^{\bullet \bullet})$ in $\mathcal B$ by: 
\[
\Phi (X^{\bullet \bullet})^{pq}:=\Phi (X^{-p,-q}),\  
d^{\prime pq}_{\Phi (X)}:=(-1)^{p+1}\Phi (d^{\prime -p-1,-q}_X),\  
d^{\prime \prime pq}_{\Phi (X)}:=(-1)^{q+1}\Phi (d^{\prime \prime -p,-q-1}_X). 
\]
If we denote $\Phi (X^{\bullet \bullet})$ by $Y^{\bullet \bullet}$ 
then $Y^{p,\bullet}=\Phi (X^{-p,\bullet})$ and $d^{\prime p,\bullet}_Y=
(-1)^{p+1}\Phi (d^{\prime -p-1,\bullet}_X)$. 

\vskip3mm
{\bf 2.6. Lemma.}\quad 
(a) \textit{If} $X^{\bullet \bullet}$ \textit{and} $Y^{\bullet \bullet}$ 
\textit{are two double complexes with} $X^{pq} = Y^{pq}$, 
$\forall p,\, q$, \textit{but with} $d_Y^{\prime} = (-1)^ad_X^{\prime}$ 
\textit{and} $d_Y^{\prime \prime} = (-1)^bd_X^{\prime \prime}$, \textit{for 
some} $a,\, b\in {\mathbb Z}$, \textit{then} $X^{\bullet \bullet} \simeq  
Y^{\bullet \bullet}$.\\  
\hspace*{3mm} (b) 
\textit{Using the notations from the last part of Definition} 2.5., 
\textit{assume that}, $\forall m\in {\mathbb Z}$, \textit{the set} 
$\{(p,q)\  \vert \  p+q=m,\  X^{pq}\neq 0\}$ \textit{is finite. Then} 
$\text{tot}(\Phi (X^{\bullet \bullet}))\simeq \Phi (\text{tot}
(X^{\bullet \bullet}))$. 
\vskip3mm

\begin{proof}
(a) $((-1)^{ap+bq}\, \text{id}_{X^{pq}})_{p,q\in {\mathbb Z}}$ is an isomorphism 
of double complexes $X^{\bullet \bullet}\overset{\sim}{\rightarrow} 
Y^{\bullet \bullet}$. 

(b) One can easily check that 
$\Phi (\text{tot}(X^{\bullet \bullet}))=\text{tot}(Z^{\bullet \bullet})$, 
where the double complex $Z^{\bullet \bullet}$ is defined by: 
\[
Z^{pq}:=\Phi (X^{-p-q}),\ d^{\prime pq}_Z:=(-1)^{p+q+1}\Phi 
(d^{\prime -p-1,-q}_X),\  d^{\prime \prime pq}_Z:=(-1)^{p+q+1}\Phi 
(d^{\prime \prime -p,-q-1}_X). 
\]
But $((-1)^{pq}\, \text{id}_{\Phi (X^{-p,-q})})_{p,q\in {\mathbb Z}}$ is an 
isomorphism of double complexes $\Phi (X^{\bullet \bullet})
\overset{\sim}{\rightarrow} Z^{\bullet \bullet}$. 
\end{proof}
    
\vskip3mm
{\bf 2.7. Lemma.}\quad 
\textit{Let} $M^{\bullet}\in \text{Ob}\, \text{C}^b(S\text{-Mod})$, 
$N^{\bullet} \in \text{Ob}\, \text{C}^b(\Lambda \text{-mod})$ 
\textit{and} $a\in {\mathbb Z}$. \textit{Then}:\\
\hspace*{3mm} (a) $\text{F}(N^{\bullet}(a))=\text{T}^a\text{F}
(N^{\bullet})(-a)$ \textit{and} $\text{G}(M^{\bullet}(a))=\text{T}^a\text{G}
(M^{\bullet})(-a)$ (\textit{one has equality, not only an isomorphism}!),\\
\hspace*{3mm} (b) $\text{F}(N^{\bullet \ast})\simeq \text{F}
(N^{\bullet})^{\vee}$,\\ 
\hspace*{3mm} (c) $\text{G}(M^{\bullet})^{\ast}\simeq 
\text{G}(M^{\bullet \ast})(-n-1)\otimes_k\overset{n+1}{\wedge}V\simeq 
\text{T}^{-n-1}\text{G}((M^{\bullet}\otimes_S{\omega}_S)^{\ast})$. 
\vskip3mm

\begin{proof}
(a) One checks, firstly, that if $M\in \text{Ob}(S\text{-Mod})$ and $N\in 
\text{Ob}(\Lambda \text{-mod})$ then $\text{F}(N)=\text{T}^a\text{F}(N)(-a)$ 
and $\text{G}(M)=\text{T}^a\text{G}(M)(-a)$. For the general case, one takes 
into account the sign convention at the end of Definition 2.1.(i).

(b) If $N\in \text{Ob}(\Lambda \text{-mod})$ then one checks easily that 
$\text{F}(N^{\ast})=\text{F}(N)^{\vee}$. Now, if $N^{\bullet}\in \text{Ob}\, 
\text{C}^b(\Lambda \text{-mod})$ then, by defintion, $\text{F}(N^{\bullet})=
\text{tot}(X^{\bullet \bullet})$ with $X^{p,\bullet}=\text{F}(N^p)$, 
$\forall p\in {\mathbb Z}$. One deduces that, using the last part of 
Definition 2.5., $\text{F}(N^{\bullet \ast})=\text{tot}
((X^{\bullet \bullet})^{\vee})$. But, by Lemma 2.6.(b), $\text{tot}
((X^{\bullet \bullet})^{\vee})\simeq \text{tot}(X^{\bullet \bullet})^{\vee}$. 

(c) If $N\in \text{Ob}(\Lambda \text{-mod})$, one can define a functor 
$\text{G}_N : S\text{-Mod} \rightarrow \text{C}(\Lambda \text{-mod})$ by 
$\text{G}_N(M)^p := M_p\otimes_kN(p)$, $d_{\text{G}_N(M)} := \sum_{i=0}^n
(X_i\cdot -)_M\otimes (e_i\cdot -)_N$. As in Definition 2.4., $\text{G}_N$ 
can be extended to a functor $\text{G}_N : \text{C}^b(S\text{-Mod}) 
\rightarrow \text{C}(\Lambda \text{-mod})$. 

Firstly, if $M\in \text{Ob}(S\text{-Mod})$ then, by Definition 2.1.(iii), 
$\text{G}_{N^{\ast}}(M^{\ast})^p \overset{\sim}{\rightarrow} 
(\text{G}_N(M)^{-p})^{\ast}$, $\forall p\in {\mathbb Z}$, and under these 
identifications, $d^p_{\text{G}_{N^{\ast}}(M^{\ast})}$ is identified with 
$(-1)^p(d^{-p-1}_{\text{G}_N(M)})^{\ast}$. Recalling the Definition 2.5., it 
follows that $\text{G}_{N^{\ast}}(M^{\ast})$ is isomorphic to a complex 
whose terms coincide with the terms of $\text{G}_N(M)^{\ast}$ but whose 
differential equals ``$-$ the differential of $\text{G}_N(M)^{\ast}$''. 
Using Lemma 2.6., one deduces now, for every complex $M^{\bullet}\in 
\text{Ob}\, \text{C}^b(S\text{-Mod})$, an isomorphism 
$\text{G}_{N^{\ast}}(M^{\bullet \ast}) \overset{\sim}{\rightarrow} 
\text{G}_N(M^{\bullet})^{\ast}$. 

Secondly, if $M\in \text{Ob}(S\text{-Mod})$ and $a\in {\mathbb Z}$ then, 
taking into account the sign convention at the end of Definition 2.1.(i), 
$\text{G}_{N(a)}(M) = \text{G}_N(M)(a)$. Using Lemma 2.6.(a) one deduces, 
for every complex $M^{\bullet}\in \text{Ob}\, \text{C}^b(S\text{-Mod})$, an 
isomorphism $\text{G}_{N(a)}(M^{\bullet}) \simeq \text{G}_N(M^{\bullet})(a)$. 

Since $\text{G} := \text{G}_{\bigwedge (V^{\ast})}$ one gets, recalling the 
isomorphisms at the end of Definition 2.1.(iv), the first isomorphism 
from the statement. The second isomorphism follows from (a).      
\end{proof} 

\vskip3mm
{\bf 2.8. Definition.} (The linear part of a minimal complex)\\
\hspace*{3mm} (i) Let $L^{\bullet}\in \text{Ob}\, \text{C}(\mathcal P)$. One 
may write 
$L^i=\bigoplus_{j\in {\mathbb Z}}S(i-j)^{\, \textstyle b_{ij}}$. For $m\in 
{\mathbb Z}$ one puts: 
\[
F_mL^i:=\underset{j\leq m}{\textstyle \bigoplus}S(i-j)^{\, \textstyle b_{ij}}.
\]
Alternatively, $F_mL^i$ is the $S$-submodule of $L^i$ generated by the 
homogeneous elements of degree $\leq m-i$. The complex $L^{\bullet}$ is 
called {\it minimal} if $\text{Im}\, d_L\subseteq S_+\cdot L^{\bullet}$. This 
is equivalent to the fact, $\forall m\in {\mathbb Z}$, $F_mL^{\bullet}:=
(F_mL^i)_{i\in {\mathbb Z}}$ is a {\it subcomplex} of $L^{\bullet}$. In this 
case, $\text{gr}_F(L^{\bullet})$ is called the {\it linear part} of 
$L^{\bullet}$. 

(ii) Similarly, let $I^{\bullet}\in \text{Ob}\, \text{C}(\mathcal I)$. One may 
write (by the last part of Definition 2.1.(iv)) $I^i=\bigoplus_{j\in 
{\mathbb Z}}\bigwedge (V^{\ast})(i-j)^{\, \textstyle c_{ij}}$. 
For $m\in {\mathbb Z}$, one 
puts:
\[
F_mI^i:=\underset{j\leq m}{\textstyle \bigoplus} {\textstyle \bigwedge}
(V^{\ast})(i-j)^{\, \textstyle c_{ij}}.
\]
Alternatively, $F_mI^i$ is the $\Lambda$-submodule of $I^i$ generated by the 
homogeneous elements of degree $\leq m-i-n-1$. The complex $I^{\bullet}$ is 
called {\it minimal} if $\text{Im}\, d_I\subseteq {\Lambda}_+\cdot 
I^{\bullet}$ or, equivalently, if $F_mI^{\bullet}:=(F_mI^i)_{i\in 
{\mathbb Z}}$ is a subcomplex of $I^{\bullet}$, $\forall m\in {\mathbb Z}$. 
In this case, $\text{gr}_F(I^{\bullet})$ is called the {\it linear part} 
of $I^{\bullet}$. 

(iii) {\it If two minimal complexes from} 
$\text{C}(\mathcal I)$ {\it are isomorphic in} 
$\text{K}(\mathcal I)$ ({\it i.e.}, {\it are homotopically equivalent}) 
{\it then they are isomorphic in} $\text{C}(\mathcal I)$ (see, for 
example, \cite{ctr}, (4.2.)). 
\vskip3mm

The following result, which is a direct consequence of Lemma A.7. from 
Appendix A, is one of the key points of the paper of 
Eisenbud, Fl\o ystad, and Schreyer \cite{efs}. 

\vskip3mm
{\bf 2.9. Lemma.}\quad (a) \textit{If} $N^{\bullet}\in \text{Ob}\, \text{C}^b
(\Lambda \text{-mod})$ \textit{then} $\text{F}(N^{\bullet})\in \text{Ob}\, 
\text{C}^b(\mathcal P)$ \textit{can be contracted to a minimal complex} 
$L^{\bullet}$ \textit{whose linear part is} $\text{F}(\text{H}^{\bullet}
(N^{\bullet}))$, \textit{where} $\text{H}^{\bullet}(N^{\bullet})$ 
\textit{is the complex with terms} $\text{H}^p(N^{\bullet})$, $p\in 
{\mathbb Z}$, \textit{and with the differential equal to} $0$. 
\textit{Moreover, this contraction induces}, 
$\forall m\in {\mathbb Z}$, \textit{a contraction of} $\text{F}({\tau}^
{\leq m}N^{\bullet})$ \textit{onto} $F_mL^{\bullet}$ \textit{and of} 
$\text{F}({\tau}^{>m}N^{\bullet})$ \textit{onto} $L^{\bullet}/F_mL^{\bullet}$. 
 
(b) \textit{If} $M^{\bullet}\in \text{Ob}\, \text{C}^b(S\text{-Mod})$ 
\textit{then} $\text{G}(M^{\bullet})\in \text{Ob}\, \text{C}(\mathcal I)$ 
\textit{can be contracted to a minimal complex} $I^{\bullet}$ 
\textit{whose linear part is} $\text{G}(\text{H}^{\bullet}(M^{\bullet}))$. 
\textit{Moreover, this contraction induces}, $\forall m\in {\mathbb Z}$, 
\textit{a contraction of} $\text{G}({\tau}^{\leq m}M^{\bullet})$ 
\textit{onto} $F_mI^{\bullet}$ \textit{and of} $\text{G}({\tau}^{>m}
M^{\bullet})$ \textit{onto} $I^{\bullet}/F_mI^{\bullet}$.
\vskip3mm
   
The next theorem is the Bernstein-Gel'fand-Gel'fand correspondence for graded 
modules. We include here a direct proof, which does not use Koszul duality. 
We use, instead, the Comparison Lemma B.1. from Appendix B. 
 
\vskip3mm
{\bf 2.10. Theorem.} (\cite{bgg}, Theorem 3.)\\
\hspace*{3mm} \textit{The functor} $\text{F} : \text{C}^b(\Lambda 
\text{-mod})\rightarrow \text{C}^b(\mathcal P)$ 
\textit{extends to an equivalence of triangulated categories} 
$\text{F} : \text{D}^b(\Lambda \text{-mod})\rightarrow \text{K}^b
(\mathcal P)$.
\vskip3mm

\begin{proof} 
If $\phi : N^{\prime \bullet}\rightarrow N^{\bullet}$ is a quasi-isomorphism 
in $\text{C}^b(\Lambda \text{-mod})$ then $\text{Con}(\phi)$ is acyclic. 
By Lemma 2.9.(a), the complex $\text{Con}(\text{F}(\phi))=\text{F}
(\text{Con}(\phi))$ is homotopically equivalent to 0, whence $\text{F}(\phi)$ 
is a homotopic equivalence. One deduces that F extends to a functor 
$\text{F} : \text{D}^b(\Lambda \text{-mod})\rightarrow \text{K}^b
(\mathcal P)$. 

We show, firstly, that this functor is {\it fully faithful}, i.e., that if 
$N^{\bullet},N^{\prime \bullet}\in \text{Ob}\, \text{C}^b(\Lambda 
\text{-mod})$ then: 
\[
\text{Hom}_{\text{D}^b(\Lambda)}(N^{\prime \bullet},N^{\bullet})
\overset{\sim}{\longrightarrow} \text{Hom}_{\text{K}^b(\mathcal P)}
(\text{F}(N^{\prime \bullet}),\text{F}(N^{\bullet})).
\tag{*}
\]
We endow $N^{\bullet}$ and $N^{\prime \bullet}$ with the filtrations 
$F^iN^{\bullet} = {\sigma}^{\geq i}N^{\bullet}$, $F^iN^{\prime \bullet} = 
{\sigma}^{\geq i}N^{\prime \bullet}$ with successive quotients 
$\text{T}^{-i}N^i$ and $\text{T}^{-i}N^{\prime i}$, respectively. If, 
$\forall i, j\in {\mathbb Z}$, the map:
\[
\text{Hom}_{\text{D}^b(\Lambda)}(N^{\prime i},\text{T}^pN^j)\longrightarrow 
\text{Hom}_{\text{K}^b(\mathcal P)}(\text{F}(N^{\prime i}),
\text{F}(\text{T}^pN^j))
\]
would be bijective for $i-j-1\leq p\leq i-j+1$ then Lemma B.1., applied to the 
functor $\text{F} : \text{D}^b(\Lambda \text{-mod})\rightarrow 
\text{K}^b(\mathcal P)$, would imply that (*) is bijective. 
Now, if $K$ and $K^{\prime}$ are two objects of $\Lambda$-mod then:
\[
\text{Hom}_{\text{D}^b(\Lambda)}(K^{\prime},K)\overset{\sim}{\longrightarrow} 
\text{Hom}_{\text{K}^b(\mathcal P)}(\text{F}(K^{\prime}),\text{F}(K))
\]
as one can easily see using the fact that $\text{Hom}_{\text{D}^b(\Lambda)}
(K^{\prime},K)\simeq \text{Hom}_{\Lambda \text{-mod}}(K^{\prime},K)$. 
Moreover, $\text{Hom}_{\text{D}^b(\Lambda)}(K^{\prime},\text{T}^pK) = 0$ for 
$p < 0$ (see (B.3.)) and, for $p > 0$, $\text{Hom}_{\text{D}^b(\Lambda)}
(K^{\prime},\text{T}^pK) = 0$ if $K$ is a direct sum of 
$\Lambda$-modules of the form 
$\bigwedge (V^{\ast})(a)$ because, in this case, $K$ is an injective object of 
$\Lambda$-mod. 

On the other hand, $\text{Hom}_{\text{K}^b(\mathcal P)}(\text{F}(K^{\prime}),
\text{F}(\text{T}^pK)) = 0$ for $p < 0$ because $\text{F}(K^{\prime})^i 
= S(i)\otimes_kK^{\prime}_i$ and $\text{F}(\text{T}^pK)^i = \text{F}(K)^{i+p} 
= S(i+p)\otimes_kK_{i+p}$, $\forall i\in {\mathbb Z}$. Moreover, if 
$L^{\bullet}\in \text{Ob}\, \text{C}^b(\mathcal P)$ and 
$\text{H}^i(L^{\bullet})_{-i} = 0$, $\forall i\in {\mathbb Z}$, then 
$\text{Hom}_{\text{K}^b(\mathcal P)}(\text{F}(K^{\prime}),L^{\bullet}) = 0$ 
because $\text{Hom}_{\text{K}^b(\mathcal P)}(\text{T}^{-i}\text{F}
(K^{\prime})^i,L^{\bullet}) \simeq \text{H}^i(L^{\bullet})_{-i}\otimes_k
(K^{\prime}_i)^{\ast} = 0$, $\forall i\in {\mathbb Z}$, and 
$\text{F}(K^{\prime})$ can be endowed with the filtration ${\sigma}^{\geq i}
\text{F}(K^{\prime})$, $i\in {\mathbb Z}$, with successive quotients 
$\text{T}^{-i}\text{F}(K^{\prime})^i$. For $L^{\bullet} = \text{F}
(\text{T}^pK)$, the condition $\text{H}^i(L^{\bullet})_{-i} = 0$, $\forall 
i\in {\mathbb Z}$, is fulfilled if $p > 0$ and $K$ is a direct sum of 
$\Lambda$-modules of the form $\bigwedge (V^{\ast})(a)$, because 
$\text{F}(\text{T}^p\bigwedge (V^{\ast})(a)) = \text{T}^{p+a}\text{F}
(\bigwedge (V^{\ast}))(-a)$ and $\text{F}(\bigwedge (V^{\ast}))$ is the 
Koszul resolution of $S/S_+$. 

Summing up, if $N^j$ is a direct sum of $\Lambda$-modules of the form 
$\bigwedge (V^{\ast})(a)$, $\forall j \leq \text{sup}\{i\in {\mathbb Z}\  
\vert \  N^{\prime i}\neq 0\}$, then Lemma B.1. implies that (*) is bijective. 
If $N^{\bullet}$ is arbitrary, one constructs, using Lemma 2.2., a 
quasi-isomorphism $N^{\bullet}\rightarrow I^{\bullet}$ with $I^{\bullet}\in 
\text{Ob}\, \text{C}^+(\mathcal I)$. For $m > \text{sup}\{i\in {\mathbb Z}\  
\vert \  N^{\prime i}\neq 0\}$ large enough, one gets a quasi-isomorphism 
$N^{\bullet}\rightarrow {\tau}^{\leq m}I^{\bullet}$. By what have been proved, 
(*) is bijective for the pair $(N^{\prime \bullet},
{\tau}^{\leq m}I^{\bullet})$, 
hence also for the pair $(N^{\prime \bullet},N^{\bullet})$.
 
Finally, the {\it essential surjectivity} can be checked as follows. By 
what have been proved, the image of $\text{F} : \text{D}^b(\Lambda 
\text{-mod})\rightarrow \text{K}^b(\mathcal P)$ is a {\it full} subcategory 
of $\text{K}^b(\mathcal P)$, closed under the functors T and $\text{T}^{-1}$ 
and under mapping cones. Moreover, $\text{F}(\text{T}^a\underline{k}(-a))=
S(a)$, $\forall a\in {\mathbb Z}$. If $L^{\bullet}\in \text{Ob}\, \text{K}^b
(\mathcal P)$ one deduces easily, by induction on $\sum_{i\in {\mathbb Z}}
\text{rk}L^i$, that $L^{\bullet}$ is isomorphic in $\text{K}^b(\mathcal P)$ to 
a complex in the image of F.
\end{proof} 
\vskip3mm

Actually, the authors of \cite{bgg} prove something more, namely that one can 
get a {\it quasi-inverse} to F by applying G and then taking convenient 
truncations (see Beilinson et al. \cite{bgs}, (2.12.) for a detailed proof). 
We shall only need the easy half of this fact, which is the 
content of the following: 

\vskip3mm
{\bf 2.11. Proposition.}\quad 
\textit{There exists a functorial quasi-isomorphism} $N^{\bullet}\rightarrow 
\text{G}\text{F}(N^{\bullet})$, $\forall N^{\bullet}\in \text{Ob}\, \text{C}^b
(\Lambda \text{-mod})$.
\vskip3mm 

\begin{proof}
We consider, firstly, the case of an object $N$ of $\Lambda$-mod. In this case 
it turns out that $\text{G}\text{F}(N)$ is an injective resolution of $N$ in 
$\Lambda$-mod. Indeed, by definition, $\text{G}\text{F}(N)=\text{tot}
(Y^{\bullet \bullet})$ with $Y^{p,\bullet}=\text{G}(\text{F}(N)^p)=\text{G}
(S(p)\otimes_kN_p)$, i.e., with $Y^{pq}=S^{p+q}(V^{\ast})\otimes_kN_p
\otimes_k\bigwedge (V^{\ast})(q)$. In particular, $\text{G}\text{F}(N)^m=0$ 
for $m<0$ and $\text{G}\text{F}(N)^0=\bigoplus_{p\in {\mathbb Z}}N_p
\otimes_k\bigwedge (V^{\ast})(-p)$. 

Let ${\beta}^p : N\rightarrow N_p\otimes_k\bigwedge (V^{\ast})(-p)$ be the 
morphism corresponding, according to Lemma 2.2., to $(-1)^p\text{id}_{N_p}$, 
and let $\beta : N\rightarrow \text{G}\text{F}(N)^0$ be the morphism 
defined by ${\beta}^p$, $p\in {\mathbb Z}$. We want to check that 
$d^0_{\text{G}\text{F}(N)}\circ \beta =0$. This is equivalent to the fact 
that, $\forall p\in {\mathbb Z}$, the diagram: 
\[
\begin{CD}
N_p\otimes_k{\textstyle \bigwedge}(V^{\ast})(-p) 
@>{\sum}X_i\otimes (e_i\cdot -)_N\otimes \text{id}>> 
V^{\ast}\otimes_kN_{p+1}\otimes_k{\textstyle \bigwedge}(V^{\ast})(-p)\\
@A{\beta}^pAA    @AA(-1)^{p+1}{\sum}X_i\otimes \text{id} \otimes 
(e_i\cdot -)_{\bigwedge (V^{\ast})}A\\
N  @>\quad \quad{\beta}^{p+1}\quad \quad>>   
N_{p+1}\otimes_k{\textstyle \bigwedge}(V^{\ast})(-p-1) 
\end{CD}
\]
anticommutes. 
According to Lemma 2.2., this is equivalent to the fact that the 
diagram: 
\[
\begin{CD}
N_p @>{\sum}X_i\otimes (e_i\cdot -)_N>> V^{\ast}\otimes N_{p+1}\\
@A(-1)^p\text{id}_{N_p}AA    @AA(-1)^{p+1}{\sum}X_i\otimes \text{id}_{N_{p+1}} 
\otimes (e_i\cdot -)_{\bigwedge (V^{\ast})}A\\
N_p @>\  \quad {\beta}^{p+1}_p\quad \  >> N_{p+1}\otimes V^{\ast}
\end{CD}
\]
anticommutes. But, according to Remark 2.3., ${\beta}^{p+1}_p = (-1)^{p+1}
\cdot (-1)^{p+1}{\sum}(e_i\cdot -)_N\otimes X_i = {\sum}(e_i\cdot -)_N
\otimes X_i$.

We have thus defined a morphism of complexes $\beta : N\rightarrow \text{G}
\text{F}(N)$. In order to show that it is a quasi-isomorphism, one may 
assume, by induction on $\text{dim}_kN$, that $N = \underline{k}(a)$ for 
some $a\in {\mathbb Z}$, and then that $a=0$, i.e., that $N = \underline{k}$. 
In this case $\text{F}(\underline{k}) = S$ and the complex $\text{G}(S)$:
\[
\cdots \rightarrow 0\rightarrow {\textstyle \bigwedge}(V^{\ast})\rightarrow 
V^{\ast}\otimes_k{\textstyle \bigwedge}(V^{\ast})(1)\rightarrow \cdots 
\rightarrow S^p(V^{\ast})\otimes_k{\textstyle \bigwedge}(V^{\ast})(p)
\rightarrow \cdots
\] 
is an injective resolution of $\underline{k}$ in $\Lambda$-mod, as one can 
easily check using the fact that the Koszul complex 
\[
0\rightarrow S(-n-1)\otimes_k\overset{n+1}{\wedge}V^{\ast}\rightarrow 
\cdots \rightarrow S(-1)\otimes_kV^{\ast}\rightarrow S\rightarrow 0
\]
is a (free) resolution of $S/S_+$ in $S$-mod.  

The general case $N^{\bullet}\in \text{Ob}\, \text{C}^b(\Lambda \text{-mod})$  
can be now deduced from the following easy observation : $\text{G}\text{F}
(N^{\bullet}) = \text{tot}(Z^{\bullet \bullet})$, where $Z^{\bullet \bullet}$ 
is the double complex with $Z^{p,\bullet} = \text{G}\text{F}(N^p)$ and with 
$d^{\prime p}_Z : Z^{p,\bullet}\rightarrow Z^{p+1,\bullet}$ equal to $\text{G}
\text{F}(d^p_N)$. Indeed, by definition, $\text{F}(N^{\bullet}) = \text{tot}
(X^{\bullet \bullet})$ with $X^{p,\bullet} = \text{F}(N^p)$ and $\text{G}
\text{F}(N^{\bullet}) = \text{tot}(Y^{\bullet \bullet})$ with $Y^{m,\bullet} = 
\text{G}(\text{tot}(X^{\bullet \bullet})^m) = \bigoplus_{p+q=m}\text{G}
(X^{pq})$. Consider the {\it triple complex} $W^{\bullet \bullet \bullet}$ 
defined by $W^{pq,\bullet} = \text{G}(X^{pq})$. We have: $Z^{p,\bullet} = 
\text{G}\text{F}(N^p) = \text{G}(X^{p,\bullet}) = \text{tot}
(W^{p,\bullet \bullet})$, hence: $\text{tot}(Z^{\bullet \bullet}) = 
\text{tot}(W^{\bullet \bullet \bullet}) = \text{tot}(Y^{\bullet \bullet}) = 
\text{G}\text{F}(N^{\bullet})$.     
\end{proof}    
\vskip3mm

{\bf 2.12. Corollary.}\quad 
$\forall N^{\bullet}\in \text{Ob}\, \text{C}^b(\Lambda\text{-mod})$, 
\textit{there exists a functorial quasi-isomorphism}: 
\[
\text{T}^{-n-1}\text{G}((\text{F}(N^{\bullet})^{\vee}\otimes_S{\omega}_S)^
{\ast})\longrightarrow N^{\bullet}. 
\]
\vskip3mm

\begin{proof}
By Proposition 2.11., there exists a functorial quasi-isomorphism 
$N^{\bullet \ast} \rightarrow \text{G}\text{F}(N^{\bullet \ast})$ and, by 
Lemma 2.7.(b), $\text{F}(N^{\bullet \ast}) \simeq \text{F}(N^{\bullet})^
{\vee}$. One gets a quasi-isomorphism $\text{G}(\text{F}(N^{\bullet})^
{\vee})^{\ast} \rightarrow N^{\bullet}$ and, by Lemma 2.7.(c), 
$\text{G}(\text{F}(N^{\bullet})^{\vee})^{\ast} \simeq \text{T}^{-n-1}\text{G} 
((\text{F}(N^{\bullet})^{\vee}\otimes_S{\omega}_S)^{\ast})$.
\end{proof}
\vskip3mm

\section{The Horrocks correspondence}

{\bf 3.1. Definition.}\quad 
(i) If $M, M^{\prime}\in \text{Ob}(S\text{-mod})$ let $I_{\mathcal P}
(M^{\prime},M)$ denote the subgroup of $\text{Hom}_{S\text{-mod}}(M^{\prime}, 
M)$ consisting of the morphisms factorizing through an object of $\mathcal P$. 
The {\it stable category} $S\text{-}\underline{\text{mod}}$ has, by 
definition, the same objects as $S$-mod, but the groups Hom are given by:
\[
\underline{\text{Hom}}_{S\text{-mod}}(M^{\prime},M) := 
\text{Hom}_{S\text{-mod}}(M^{\prime},M)/I_{\mathcal P}(M^{\prime},M). 
\]

(ii) Similarly, using the full subcategory $\widetilde {\mathcal P}$ of 
$\text{Coh}\, {\mathbb P}$ consisting of finite direct sums of invertible 
sheaves ${\mathcal O}_{\mathbb P}(a)$, $a\in {\mathbb Z}$, one defines the 
{\it stable category} $\underline{\text{Coh}}\, {\mathbb P}$. 
\vskip3mm 

{\bf 3.2. Definition.}\quad 
A complex $K^{\bullet}\in \text{Ob}\, \text{C}^b(\mathcal P)$ is called a 
{\it Horrocks complex} if it satisfies the following equivalent conditions: 

(1) $\text{H}^i(K^{\bullet}) = 0$ \textit{for} $i\leq -2$ \textit{and} 
$\text{H}^i(K^{\bullet \vee}) = 0$ \textit{for} $i\leq 1$\\
\hspace*{3mm} (2) $\text{H}^i(K^{\bullet \vee}) = 0$ \textit{for} $i\leq 1$ 
\textit{and} $\text{dim}\, \text{H}^i(K^{\bullet \vee})\leq n+2-i$, 
\textit{for} $i>1$\\
\hspace*{3mm} (3) $\text{H}^i(K^{\bullet}) = 0$ \textit{for} $i\leq -2$ 
\textit{and} $\text{dim}\, \text{H}^i(K^{\bullet})\leq n-1-i$, \textit{for} 
$i\geq -1$. 

Here ``dim'' stands for ``Krull dimension''. The equivalence of these 
conditions follows from Lemma 3.3. below. Condition (1) implies that if 
$K^{\prime \bullet}\in \text{Ob}\, \text{C}^b(\mathcal P)$ is homotopically 
equivalent to a Horrocks complex then it is a Horrocks complex. 

Let $M\in \text{Ob}(S\text{-mod})$ and let $L^{\bullet}$ (resp., 
$L^{\prime \bullet}$) be a free resolution of $M$ (resp., $M^{\vee}$) in 
$S$-mod. One can concatenate the complexes $L^{\prime \bullet}$ and 
$\text{T}^{-1}(L^{\bullet \vee})$ using the composite morphism 
$L^{\prime 0} \twoheadrightarrow M^{\vee} \hookrightarrow L^{0 \vee}$. 
The dual $K^{\bullet}$ of the resulting complex is a Horrocks complex.  
We call it a {\it Horrocks resolution} of $M$. 
\vskip3mm 

{\bf 3.3. Lemma.}\quad 
\textit{Let} $A$ \textit{be a Noetherian} (\textit{commutative}) \textit{ring 
and let} $P^{\bullet}$ \textit{be a left bounded complex of finitely 
generated projective} $A$-\textit{modules}. \textit{Then the following 
conditions are equivalent}:\\
\hspace*{3mm} (i) $\text{H}^i(P^{\bullet}) = 0$, $\forall i < 0$\\
\hspace*{3mm} (ii) $\forall i > 0$, $\forall \mathfrak{p}\in \text{Supp}\, 
\text{H}^i(P^{\bullet \vee})\subseteq \text{Spec}A$, $\text{depth}\, 
A_{\mathfrak{p}}\geq i$.
\vskip3mm

\begin{proof}
(i)$\Rightarrow$(ii) Let $i > 0$ and let $\mathfrak{p}\in \text{Spec}\, A$ 
with $\text{depth}\, A_{\mathfrak{p}} < i$. Let $M := C^0(P^{\bullet}) := 
\text{Coker}(P^{-1}\rightarrow P^0)$. Condition (i) implies that $M$ has 
finite projective dimension. Now, the Auslander-Buchsbaum formula implies that 
the projective dimension of the $A_{\mathfrak{p}}$-module $M_{\mathfrak{p}}$ 
is $\leq \text{depth}\, A_{\mathfrak{p}} < i$, hence $\text{H}^i
(P^{\bullet \vee})_{\mathfrak{p}} \simeq 
\text{Ext}^i_{A_{\mathfrak{p}}}(M_{\mathfrak{p}},A_{\mathfrak{p}}) = 0$, 
whence $\mathfrak{p}\notin \text{Supp}\, \text{H}^i(P^{\bullet \vee})$. 

(ii)$\Rightarrow$(i) We use induction on $m := \text{sup}\{i\in {\mathbb Z}\  
\vert\  \text{H}^i(P^{\bullet \vee})\neq 0\}$. The case $m\leq 0$ is clear. 
For the proof of the induction step $(m-1)\rightarrow m$, consider the 
$A$-module $N := C^{-1}(P^{\bullet}) := \text{Coker}(P^{-2}\rightarrow 
P^{-1})$. Applying the induction hypothesis to $\text{T}^{-1}P^{\bullet}$, one 
gets that $\text{H}^i(P^{\bullet}) = 0$, $\forall i < -1$, hence the 
sequence:
\[
0\rightarrow P^{-r}\rightarrow \cdots \rightarrow P^{-2}\rightarrow P^{-1}
\rightarrow N\rightarrow 0
\]
is exact. We {\it assert} that $\text{Ass}(N)\subseteq \text{Ass}(A)$. 
Indeed, let $\mathfrak{p}\in \text{Ass}(N)$ and let $d := \text{depth}\, 
A_{\mathfrak{p}}$. It follows from the Auslander-Buchsbaum formula that the 
projective dimension of the $A_{\mathfrak{p}}$-module $N_{\mathfrak{p}}$ is 
$d$, which implies that $\text{Ext}^d_{A_{\mathfrak{p}}}(N_{\mathfrak{p}}, 
A_{\mathfrak{p}}) \neq 0$. If $d > 0$ then $\text{Ext}^d_{A_{\mathfrak{p}}}
(N_{\mathfrak{p}},A_{\mathfrak{p}})\simeq \text{H}^{d+1}(P^{\bullet \vee})_
{\mathfrak{p}}$ hence, by (ii), $\text{depth}\, A_{\mathfrak{p}} \geq d+1$, 
a {\it contradiction}. It remains that $d = 0$, i.e., $\mathfrak{p}\in 
\text{Ass}(A)$.

Now, $\text{H}^{-1}(P^{\bullet})\simeq \text{Ker}(N\rightarrow P^0)$, hence 
$\text{Ass}(\text{H}^{-1}(P^{\bullet}))\subseteq \text{Ass}(N)\subseteq 
\text{Ass}(A)$. If $\mathfrak{p}\in \text{Ass}(A)$ then, by (ii), the 
sequence:
\[
P^{0\vee}_{\mathfrak{p}}\rightarrow P^{-1\vee}_{\mathfrak{p}}\rightarrow 
\cdots \rightarrow P^{-r\vee}_{\mathfrak{p}}\rightarrow 0
\]
is exact. Since it consists of free $A_{\mathfrak{p}}$-modules, its dual is 
also exact. In particular, it follows that $\text{H}^{-1}(P^{\bullet})_
{\mathfrak{p}} = 0$. One deduces that $\text{Ass}(\text{H}^{-1}(P^{\bullet})) 
= \emptyset$, i.e., $\text{H}^{-1}(P^{\bullet}) = 0$. 
\end{proof}
\vskip3mm 

{\bf 3.4. Theorem.}\quad 
\textit{The functor} $C^{-1} : \text{C}^b(\mathcal P)\rightarrow S\text{-mod}$ 
\textit{associating to a complex} $L^{\bullet}$ \textit{the cokernel of the 
differential} $d^{-2}_L : L^{-2}\rightarrow L^{-1}$ \textit{induces a 
functor} $\underline{C}^{-1} : \text{K}^b(\mathcal P)\rightarrow 
S\text{-}\underline{\text{mod}}$ 
\textit{which}, \textit{restricted to the full subcategory} 
$\mathcal H$ \textit{of} $\text{K}^b(\mathcal P)$ \textit{consisting of 
Horrocks complexes}, \textit{is an equivalence of categories}.
\vskip3mm

\begin{proof}
If a morphism $f : L^{\prime \bullet}\rightarrow L^{\bullet}$ in 
$\text{C}^b(\mathcal P)$ is homotopic to 0 then $C^{-1}(L^{\prime 
\bullet})\rightarrow C^{-1}(L^{\bullet})$ factorizes through $L^{-1}
\twoheadrightarrow C^{-1}(L^{\bullet})$ (and through $C^{-1}(L^{\prime 
\bullet})\rightarrow L^{\prime 0}$) hence $C^{-1}$ induces a functor 
$\underline{C}^{-1}$ as in the statement. 

The fact that $\underline{C}^{-1}\, \vert \, {\mathcal H}$ is 
{\it fully faithful} follows from the more general Lemma 3.5. below. The 
fact that $\underline{C}^{-1}\, \vert \, {\mathcal H}\rightarrow S\text{-}
\underline{\text{mod}}$ is {\it essentially surjective} was already 
observed in the last part of Definition 3.2..
\end{proof}
\vskip3mm

{\bf 3.5. Lemma.}\quad 
\textit{Let} $L^{\bullet},L^{\prime \bullet}\in \text{Ob}\, \text{C}
(\mathcal P)$. \textit{If} $\text{H}^i(L^{\bullet}) = 0$ \textit{for} 
$i\leq -2$ \textit{and if} $\text{H}^i(L^{\prime \bullet \vee}) = 0$ 
\textit{for} $i\leq 1$, \textit{then the morphism}: 
\[
\text{Hom}_{\text{C}(\mathcal P)}(L^{\prime \bullet},L^{\bullet})
\longrightarrow \text{Hom}_{S\text{-mod}}(C^{-1}(L^{\prime 
\bullet}),C^{-1}(L^{\bullet}))
\]
\textit{is surjective and induces an isomorphism}:
\[
\text{Hom}_{\text{K}(\mathcal P)}(L^{\prime \bullet},L^{\bullet})
\overset{\sim}{\longrightarrow} \underline{\text{Hom}}_{S\text{-mod}}
(C^{-1}(L^{\prime \bullet}),C^{-1}(L^{\bullet})).
\]

\begin{proof}
The complex $\cdots \rightarrow L^{-2}\rightarrow L^{-1}\rightarrow 0$ is a 
free resolution of $C^{-1}(L^{\bullet})$ in $S$-mod, and the complex 
$\cdots \rightarrow L^{\prime 1\vee}\rightarrow L^{\prime 0\vee}\rightarrow 
0$ is a free resolution of $C^{-1}(L^{\prime \bullet})^{\vee}$. 
Now, one uses the following two elementary facts: (1) if $P^{\bullet}\in 
\text{Ob}\, \text{C}^{\leq 0}(\mathcal P)$, 
$M^{\bullet}\in \text{Ob}\, \text{C}^{\leq 0}(S\text{-mod})$ 
and if $\text{H}^i(M^{\bullet}) = 0$ for $i < 0$, 
then any morphism $C^0(P^{\bullet})\rightarrow C^0(M^{\bullet})$ can be 
lifted to a morphism of complexes $P^{\bullet}\rightarrow M^{\bullet}$; (2) 
if, moreover, $C^0(P^{\bullet})\rightarrow C^0(M^{\bullet})$ factorizes 
through an object of $\mathcal P$ (hence through $M^0\twoheadrightarrow 
C^0(M^{\bullet})$) then the morphism of augmented complexes:
\[
\begin{CD}
\cdots @>>> P^{-1} @>>> P^0 @>>> C^0(P^{\bullet}) @>>> 0\\
           @.  @VVV        @VVV     @VVV\\
\cdots @>>> M^{-1} @>>> M^0 @>>> C^0(M^{\bullet}) @>>> 0
\end{CD}
\]
is homotopic to 0.
\end{proof}

\vskip3mm
{\bf 3.6. Theorem.}\quad 
\textit{For} $N^{\bullet}\in \text{Ob}\, \text{C}^b(\Lambda \text{-mod})$ 
\textit{the complex} $\text{F}(N^{\bullet})$ \textit{is a Horrocks complex 
if and only if the linear part of a minimal free resolution of} $N^{\bullet}$ 
\textit{in} $\Lambda$-mod \textit{is of the form} $\bigoplus_{i=-1}^{n-1}
\text{T}^{-i}\text{G}(H^i)$, \textit{where} $H^i$ \textit{is the} 
$k$-\textit{vector space graded dual of a finitely generated graded} 
$S$-\textit{module of Krull dimension} $\leq i+1$, $i = -1,\ldots ,n-1$. 
\vskip3mm

\begin{proof}
Let us denote $\text{F}(N^{\bullet})$ by $K^{\bullet}$. By Corollary 2.12. and 
by Lemma 2.9.(b), the linear part of a minimal free resolution of 
$N^{\bullet}$ is isomorphic to $\bigoplus_{i\in {\mathbb Z}}\text{T}^{-i}
\text{G}(H^i)$, where $H^i = \text{H}^i(\text{T}^{-n-1}
((K^{\bullet \vee}\otimes_S{\omega}_S)^{\ast})) \simeq 
(\text{H}^{n+1-i}(K^{\bullet \vee}\otimes_S{\omega}_S))^{\ast}$. 
One can now conclude, using condition (2) from Definition 3.2.. 
\end{proof}

\vskip3mm
{\bf 3.7. Definition.}\quad 
A {\it minimal} complex $G^{\bullet}\in \text{Ob}\, \text{C}^-(\mathcal I)$ 
with $\text{H}^p(G^{\bullet}) = 0$ for $p << 0$ is called a 
{\it Horrocks-Trautmann complex} if it satisfies the following conditions 
(compare with \cite{ctr}, (1.6.)): 

(1) $F_{n-1}G^{\bullet} = G^{\bullet}$ and $F_0G^{\bullet} = 0$, i.e., 
$G^p \simeq {\bigoplus}_{i=1}^{n-1}\bigwedge 
(V^{\ast})(p-i)^{\textstyle c_{pi}}$, 
$\forall p\in {\mathbb Z}$,\\
\hspace*{3mm} (2) $\underset{p\rightarrow \infty}{\text{lim}}
(c_{-p,i}/p^{i+1}) = 0$, $i = 1,\dots ,n-1$. 
\vskip3mm 

{\bf 3.8. Lemma.}\quad 
\textit{A minimal complex} $G^{\bullet}\in \text{Ob}\, \text{C}^-(\mathcal I)$ 
\textit{is a Horrocks-Trautmann complex if and only if its linear part is of 
the form} $\bigoplus_{i=1}^{n-1}\text{T}^{-i}\text{G}(H^i)$, \textit{where} 
$H^i$ \textit{is the} $k$-\textit{vector space graded 
dual of a finitely generated graded} 
$S$-\textit{module of Krull dimension} $\leq i+1$, $i = 1,\dots ,n-1$.
\vskip3mm

\begin{proof}
The equivalence can be proved by applying to $G^{\bullet \ast}$ the following: 

{\bf Assertion.}\quad  
\textit{For a minimal complex} $I^{\bullet}\in \text{Ob}\, 
\text{C}^+(\mathcal I)$, \textit{the following conditions are equivalent}:\\
\hspace*{3mm} (i) $\text{H}^p(I^{\bullet}) = 0$, \textit{for} $p >> 0$\\
\hspace*{3mm} (ii) \textit{The linear part of} $I^{\bullet}$ \textit{is of the 
form} $\bigoplus_{i\in {\mathbb Z}}\text{T}^{-i}\text{G}(M^i)$, 
\textit{where} $M^i$ \textit{is a finitely generated graded} 
$S$-\textit{module}, $\forall i\in {\mathbb Z}$.

(i)$\Rightarrow$(ii) $\text{H}^p(I^{\bullet}) = 0$, $\forall p > m$, for some 
$m\in {\mathbb Z}$. Let $Z^m := \text{Ker}(I^m\rightarrow I^{m+1})$.  
${\sigma}^{\geq m}I^{\bullet}$ is a minimal right free resolution 
of $\text{T}^{-m}Z^m$ in $\Lambda$-mod. 
By Proposition 2.11., $\text{G}\text{F}(\text{T}^{-m}Z^m)$ is a 
right free resolution of $\text{T}^{-m}Z^m$ in $\Lambda$-mod.  
From Lemma 2.9.(b), it can be contracted to a minimal complex $J^{\bullet}$ 
in $\text{C}^+(\mathcal I)$, whose linear part is 
$\bigoplus_{i\in {\mathbb Z}}\text{T}^{-i}\text{G}(\text{H}^i(\text{F}
(\text{T}^{-m}Z^m))) = \bigoplus_{i\in {\mathbb Z}}\text{T}^{-i}\text{G}
(\text{H}^{i-m}(\text{F}(Z^m)))$. ${\sigma}^{\geq m}I^{\bullet}$ and 
$J^{\bullet}$ are isomorphic in $\text{D}^+(\Lambda \text{-mod})$ hence, 
since every free object of $\Lambda$-mod is an injective object of this 
category, they are isomorphic in $\text{K}^+(\mathcal I)$ and consequently, 
by (2.8.)(iii), isomorphic in $\text{C}^+(\mathcal I)$. One deduces that the 
linear part of ${\sigma}^{\geq m}I^{\bullet}$ is isomorphic to 
$\bigoplus_{i\in {\mathbb Z}}\text{T}^{-i}\text{G}(\text{H}^{i-m}(\text{F}
(Z^m)))$.   

(ii)$\Rightarrow$(i) It suffices to prove that if $M\in \text{Ob}
(S\text{-mod})$ then $\text{H}^p(\text{G}(M)) = 0$ for $p >> 0$. Let 
$L^{\bullet}$ be a finite free resolution of $M$ in $S$-mod. $\text{G}(M)$ and 
$\text{G}(L^{\bullet})$ are quasi-isomorphic (even homotopically equivalent). 
Since $\text{G}(S)$ is a right free resolution of $\underline{k}$ in 
$\Lambda$-mod (see the proof of Proposition 2.11.) it follows that 
$\text{H}^p(\text{G}(L^{\bullet})) = 0$ for $p >> 0$.
\end{proof}
\vskip3mm       

{\bf 3.9. Theorem.}\quad 
\textit{There exists an equivalence of categories between the full 
subcategory} $\mathcal{HT}$ 
\textit{of} $\text{K}^-(\mathcal I)$ 
\textit{consisting of Horrocks-Trautmann complexes and the full subcategory 
of} $\underline{\text{Coh}}\, {\mathbb P}(V)$ 
\textit{consisting of the coherent sheaves} 
$\mathcal F$ \textit{with the property that} $\text{H}^0
({\mathcal F}(-t)) = 0$, \textit{for} $t >> 0$. 
\vskip3mm

\begin{proof}
The equivalence from the statement will appear as a composition of 
previously established equivalences. 

(1) Let $\mathcal B$ be the full subcategory of $\text{Coh}\, {\mathbb P}(V)$ 
consisting of the coherent sheaves with the property from the statement. Let 
$\mathcal A$ be the full subcategory of $S$-mod consisting of the modules 
of projective dimension $\leq n-1$. Using Theorem 1.1. (Graded Serre Duality) 
one sees that the functor $(-)\sptilde : S\text{-mod}\rightarrow \text{Coh}\, 
{\mathbb P}(V)$ induces an equivalences of categories between $\mathcal A$ 
and $\mathcal B$. Moreover, this equivalence induces an equivalence between 
the correponding full subcategories $\underline{\mathcal A}$ and 
$\underline{\mathcal B}$ of $S\text{-}\underline{\text{mod}}$ and 
$\underline{\text{Coh}}\, {\mathbb P}(V)$, respectively. 

(2) By Lemma 1.2., 
the equivalence $\underline{C}^{-1} : {\mathcal H}\rightarrow S\text{-}
\underline{\text{mod}}$ from Theorem 3.4. induces an equivalence between the 
full subcategory ${\mathcal H}^{\prime}$ of $\mathcal H$ consisting of the 
Horrocks complexes $K^{\bullet}$ with the additional property that 
$\text{H}^i(K^{\bullet \vee}) = 0$ for $i\geq n+1$ and 
$\underline{\mathcal A}$. 

(3) Finally, there is a well-known equivalence $\Phi$ between the full 
subcategory $\mathcal K$ of $\text{K}^-(\mathcal I)$ consisting of the 
complexes $I^{\bullet}$ with $\text{H}^p(I^{\bullet}) = 0$ for $p << 0$ and 
$\text{D}^b(\Lambda \text{-mod})$. $\Phi$ associates to $I^{\bullet}$ a 
convenient truncation ${\tau}^{\geq m}I^{\bullet}$ with $m << 0$ (it 
suffices that $\text{H}^p(I^{\bullet}) = 0$ for $p < m$) and its 
quasi-inverse associates to a complex in $\text{D}^b(\Lambda \text{-mod})$ a 
free resolution of it. Now, by Theorem 3.6. and Lemma 3.8., 
the composition of the BGG equivalence (Theorem 2.10.)  
$\text{F} : \text{D}^b(\Lambda \text{-mod})\rightarrow \text{K}^b(\mathcal P)$ 
and $\Phi$ induces an equivalence between $\mathcal{HT}$ and 
${\mathcal H}^{\prime}$. 
\end{proof}
\vskip3mm    

{\bf 3.10. Example.} (Eilenberg-MacLane sheaves)\\
\hspace*{3mm} Let $0 < i < n$ and let $E$ be a finitely generated graded 
$S$-module of Krull dimension $\leq i+1$. Consider a minimal free resolution 
of $E$ in $S$-mod: 
\[
0\rightarrow Q^{-n-1}\rightarrow \cdots \rightarrow Q^0\rightarrow E
\rightarrow 0.
\]
Applying Lemma 3.3. to $P^{\bullet} := \text{T}^{n-i}(Q^{\bullet \vee})$ one 
derives that $\text{H}^j(Q^{\bullet \vee}) = 0$ for $j\leq n-i-1$. Using 
condition (1) from Definition 3.2. one deduces that $K^{\bullet} := 
\text{T}^{n-i+1}(Q^{\bullet \vee})$ is a Horrocks complex. It is a Horrocks 
resolution of $M := \text{Coker}((Q^{-n+i+1})^{\vee}\rightarrow 
(Q^{-n+i})^{\vee})$. $M$ has a minimal free resolution:
\[
0\rightarrow (Q^0)^{\vee}\rightarrow \cdots \rightarrow (Q^{-n+i+1})^{\vee}
\rightarrow (Q^{-n+i})^{\vee}\rightarrow M\rightarrow 0.
\]
Let ${\mathcal F} := {\widetilde M}$. Since $\text{H}^j(K^{\bullet \vee}) 
= 0$ for $j \neq n-i+1$ and $\text{H}^{n-i+1}(K^{\bullet \vee}) \simeq E$, 
it follows from the proof of Theorem 3.6. that the Horrocks-Trautmann complex 
associated to $\mathcal F$ is $\text{T}^{-i}\text{G}(H)$ where 
$H = (E\otimes_S{\omega}_S)^{\ast}$. Moreover, by 
Graded Serre Duality (Theorem 1.1.), $\text{H}^j_{\ast}{\mathcal F} = 0$ for 
$0 < j < n$, $j\neq i$, and $\text{H}^i_{\ast}{\mathcal F} \simeq H$. 

When $E$ is of finite length $\mathcal F$ is a locally free sheaf. The locally 
free sheaves of this kind were called {\it Eilenberg-MacLane bundles} in 
Horrocks \cite{ho2}.

\section{The Horrocks correspondence and the BGG correspondence}

{\bf 4.1. Definition.} The {\it geometric BGG functor} is the functor 
$\text{L} : \Lambda\text{-mod}\rightarrow \text{C}^b(\text{Coh}\, 
{\mathbb P})$ defined by $\text{L}(N) := \text{F}(N)\sptilde$. 

We denote by $\Lambda\text{-}\underline{\text{mod}}$ the {\it stable category} 
of $\Lambda$-mod with respect to its full subcategory $\mathcal I$ consisting 
of free objects (see Definition 3.1.). 
\vskip3mm 

{\bf 4.2. Lemma.}\quad 
\textit{If} $N,N^{\prime}\in \text{Ob}(\Lambda\text{-mod})$ \textit{then}, 
$\forall p\geq 1$:
\[
\text{Hom}_{\text{K}^b(\mathbb P)}(\text{L}(N^{\prime}),\text{T}^p\text{L}(N)) 
\overset{\sim}{\longrightarrow}\text{Hom}_{\text{D}^b(\mathbb P)}
(\text{L}(N^{\prime}),\text{T}^p\text{L}(N)).
\]
\begin{proof}
The lemma is an immediate application of Lemma B.4., taking 
into account that $\text{H}^i({\mathcal O}_{\mathbb P}(a)) = 0$ for 
$i > 0$, $i\neq n$, $\forall a\in {\mathbb Z}$, and $\text{H}^n({\mathcal O}_
{\mathbb P}(a)) = 0$, $\forall a\geq -n$. 
\end{proof}
\vskip3mm

{\bf 4.3. Corollary.}\quad 
\textit{If} $N, N^{\prime}\in \text{Ob}(\Lambda\text{-mod})$ \textit{then}, 
$\forall p\geq 1$:
\[
\text{Hom}_{\text{D}^b(\Lambda)}(N^{\prime},\text{T}^pN)
\overset{\sim}{\longrightarrow}\text{Hom}_{\text{D}^b(\mathbb P)}
(\text{L}(N^{\prime}),\text{T}^p\text{L}(N)).
\]
\begin{proof}
By the BGG correspondence for graded modules Theorem 2.10.:
\[
\text{Hom}_{\text{D}^b(\Lambda)}(N^{\prime},\text{T}^pN)
\overset{\sim}{\longrightarrow}\text{Hom}_{\text{K}^b(\mathcal P)}
(\text{F}(N^{\prime}),\text{T}^p\text{F}(N))
\]
and, on the other hand, it is obvious that:
\[
\text{Hom}_{\text{K}^b(\mathcal P)}(\text{F}(N^{\prime}),\text{T}^p\text{F}
(N))\overset{\sim}{\longrightarrow}\text{Hom}_{\text{K}^b(\mathbb P)}
(\text{L}(N^{\prime}),\text{T}^p\text{L}(N)).
\]
It only remains, now, to apply Lemma 4.2..
\end{proof}
\vskip3mm

The following theorem is the Bernstein-Gel'fand-Gel'fand correspondence for 
coherent sheaves on projective spaces. We include here a direct proof of 
this result.

\vskip3mm
{\bf 4.4. Theorem.} (\cite{bgg}, Theorem 2.)\quad 
\textit{The functor} $\text{L} : \Lambda\text{-mod}\rightarrow 
\text{C}^b(\text{Coh}\, {\mathbb P})$ \textit{induces an equivalence of 
categories} $\text{L} : \Lambda\text{-}\underline{\text{mod}}\rightarrow 
\text{D}^b(\text{Coh}\, {\mathbb P})$. 
\vskip3mm

\begin{proof}
$\text{L}(\bigwedge (V^{\ast}))$ is the tautological Koszul complex on 
${\mathbb P}(V)$:
\[
0\rightarrow {\mathcal O}_{\mathbb P}(-n-1)\otimes_k\overset{n+1}{\wedge}
V^{\ast}\rightarrow \cdots \rightarrow {\mathcal O}_{\mathbb P}(-1)
\otimes_kV^{\ast}\rightarrow {\mathcal O}_{\mathbb P}\rightarrow 0
\]
hence, if $P$ is a free object of $\Lambda$-mod then $\text{L}(P)$ is an 
acyclic complex. It follows that $\text{L} : \Lambda\text{-mod}\rightarrow 
\text{C}^b(\text{Coh}\, {\mathbb P})$ induces a functor $\text{L} : \Lambda 
\text{-}\underline{\text{mod}}\rightarrow \text{D}^b(\text{Coh}\, 
{\mathbb P})$. 

We firstly show that the induced functor is {\it fully faithful}. Let 
$N, N^{\prime}\in \text{Ob}(\Lambda\text{-mod})$. We have to show that the 
morphism: 
\[
\text{Hom}_{\Lambda\text{-mod}}(N^{\prime},N)\longrightarrow 
\text{Hom}_{\text{D}^b(\mathbb P)}(\text{L}(N^{\prime}), \text{L}(N))
\tag{*}
\]
is surjective and that its kernel consists of the morphisms  factorizing 
through a free object of $\Lambda$-mod. Consider an exact sequence 
$0 \rightarrow K \rightarrow P \rightarrow N \rightarrow 0$ with $P$ a free 
object of $\Lambda$-mod. From Lemma 4.2.:
\[
\text{Hom}_{\text{K}^b(\mathbb P)}(\text{L}(N^{\prime}),\text{T}\text{L}(P)) 
\overset{\sim}{\longrightarrow}
\text{Hom}_{\text{D}^b(\mathbb P)}(\text{L}(N^{\prime}),\text{T}\text{L}(P))
\]
and $\text{Hom}_{\text{D}^b(\mathbb P)}(\text{L}(N^{\prime}),\text{T}\text{L}
(P)) = 0$ since $\text{L}(P)$ is acyclic.  
Now, applying $\text{Hom}_{\text{K}^b(\mathbb P)}
(\text{L}(N^{\prime}),-)$ and $\text{Hom}_{\text{D}^b(\mathbb P)}
(\text{L}(N^{\prime}),-)$ to the complex in $\text{K}^b(\mathbb P)$:
\[
\text{L}(P)\rightarrow \text{L}(N)\rightarrow \text{T}\text{L}(K)\rightarrow 
\text{T}\text{L}(P)
\]
deduced (see \cite{coa}, (2)(i),(ii)) from the semi-split short exact 
sequence:
\[
0\rightarrow \text{L}(K)\rightarrow \text{L}(P)\rightarrow \text{L}(N)
\rightarrow 0,
\]
one gets a commutative diagram with exact rows:
\[
\begin{CD}
\text{Hom}_{\text{K}(\mathbb P)}(\text{L}(N^{\prime}),\text{L}(P))  @>>>
\text{Hom}_{\text{K}(\mathbb P)}(\text{L}(N^{\prime}),\text{L}(N))  @>>>
\text{Hom}_{\text{K}(\mathbb P)}(\text{L}(N^{\prime}),\text{T}\text{L}(K))
{\, \rightarrow 0}\\
@VVV  @VVV  @VV{\wr}V @.\\
\text{Hom}_{\text{D}(\mathbb P)}(\text{L}(N^{\prime}),\text{L}(P))  @>>>
\text{Hom}_{\text{D}(\mathbb P)}(\text{L}(N^{\prime}),\text{L}(N))  @>>> 
\text{Hom}_{\text{D}(\mathbb P)}(\text{L}(N^{\prime}),\text{T}\text{L}(K)) 
{\, \rightarrow 0}
\end{CD}
\]
By Lemma 4.2., the vertical arrow from the right hand side of the diagram is 
an isomorphism. Moreover, $\text{Hom}_{\text{D}(\mathbb P)}(\text{L}
(N^{\prime}),\text{L}(P)) = 0$ because $\text{L}(P)$ is acyclic. Since the 
vertical arrows in the commutative diagram: 
\[
\begin{CD}
\text{Hom}_{\Lambda\text{-mod}}(N^{\prime},P) @>>>  
\text{Hom}_{\Lambda\text{-mod}}(N^{\prime},N)\\
@VV{\wr}V @VV{\wr}V\\
\text{Hom}_{\text{K}(\mathbb P)}(\text{L}(N^{\prime}),\text{L}(P))
@>>>
\text{Hom}_{\text{K}(\mathbb P)}(\text{L}(N^{\prime}),\text{L}(N))
\end{CD}
\]
are clearly isomorphisms, one deduces that the morphism (*) is surjective and 
that its kernel consists of the morphisms factorizing through $P\rightarrow 
N$.

The {\it essential surjectivity} of $\text{L} : 
\Lambda\text{-}\underline{\text{mod}}\rightarrow 
\text{D}^b(\mathbb P)$ can now be checked, in a well-known manner, using the 
following observations:

(1) By what have been proved, the image of L is a {\it full} subcategory of 
$\text{D}^b(\mathbb P)$. 

(2) If $N\in \text{Ob}(\Lambda\text{-mod})$ and one considers a short exact 
sequence $0 \rightarrow N \rightarrow I \rightarrow Q \rightarrow 0$ with 
$I$ a free object of $\Lambda\text{-mod}$ then the connecting morphism 
$w : \text{L}(Q)\rightarrow \text{T}\text{L}(N)$ deduced (see \cite{coa}, 
(2)(ii)) from the semi-split short exact sequence:
\[
0\rightarrow \text{L}(N)\rightarrow \text{L}(I)\rightarrow \text{L}(Q)
\rightarrow 0
\]
is a quasi-isomorphism because $\text{L}(I)$ is acyclic. Similarly, 
considering a short exact sequence $0 \rightarrow K \rightarrow P \rightarrow 
N \rightarrow 0$ with $P$ a free object of $\Lambda$-mod one gets a 
quasi-isomorphism $\text{T}^{-1}\text{L}(N) \rightarrow \text{L}(K)$.

(3) Let $u : N^{\prime} \rightarrow N$ be a morphism in $\Lambda$-mod. 
Consider an embedding $v : N^{\prime}\rightarrow I^{\prime}$ of 
$N^{\prime}$ into a free object $I^{\prime}$ of $\Lambda$-mod and 
define $C\in \text{Ob}(\Lambda\text{-mod})$ by the  
short exact sequence:
\[
0\rightarrow N^{\prime}\overset{(u,v)}{\longrightarrow} 
N\oplus I^{\prime}\longrightarrow C\rightarrow 0.
\]
By applying L to this short exact sequence one gets a semi-split short exact 
sequence, hence $\text{L}(C)$ is homotopically equivalent to $\text{Con}\, 
\text{L}((u,v))$. Moreover, $\text{Con}\, \text{L}((u,v))$ is 
quasi-isomorphic to $\text{Con}\, \text{L}(u)$ because $\text{L}(I^{\prime})$ 
is acyclic, whence one gets a quasi-isomorphism $\text{L}(C) \rightarrow 
\text{Con}\, \text{L}(u)$.  

(4) $\text{L}(\underline{k}(-a)) = \text{T}^{-a}{\mathcal O}_{\mathbb P}(a)$, 
$\forall a\in {\mathbb Z}$. 

Using these observations and the fact that every coherent sheaf on 
${\mathbb P}(V)$ has a finite resolution with finite direct sums of 
invertible sheaves ${\mathcal O}_{\mathbb P}(a)$, one deduces immediately 
(as, for example, in the last part of the proof of \cite{coa}, Theorem 7) 
that $\text{L} : \Lambda\text{-}\underline{\text{mod}}\rightarrow 
\text{D}^b(\mathbb P)$ is essentially surjective.  
\end{proof}
\vskip3mm

{\bf 4.5. Corollary.} (\cite{bgg}, Remark 3 after Theorem 1)\quad 
\textit{For every} ${\mathcal F}^{\bullet}\in \text{Ob}\, 
\text{C}^b(\text{Coh}\, {\mathbb P})$ \textit{there exists} 
$N\in \text{Ob}(\Lambda\text{-mod})$ 
\textit{annihilated by} $\text{soc}(\Lambda) = \overset{n+1}{\wedge}V$ 
\textit{such that} ${\mathcal F}^{\bullet} \simeq \text{L}(N)$ \textit{in} 
$\text{D}^b(\text{Coh}\, {\mathbb P})$. \textit{Moreover}, $N$ \textit{is 
unique up to isomorphism}.
\vskip3mm

For a proof see, for example, \cite{coa},(8).    

\vskip3mm
{\bf 4.6. Corollary.}\quad 
\textit{If} ${\mathcal F}^{\bullet}$ \textit{and} $N$ \textit{are as in 
Corollary} 4.5. \textit{then}, $\forall i\in {\mathbb Z}$, 
$\text{H}^i(\text{F}(N)) \simeq 
\bigoplus_{j\geq -i}{\mathbb H}^i({\mathcal F}^{\bullet}(j))$ \textit{as} 
$S$-\textit{modules} (\textit{where} 
$\mathbb H$ \textit{denotes the hypercohomology}).
\vskip3mm

\begin{proof}
$\text{H}^i(\text{F}(N))_j = 0$ for $j < -i$ because $\text{F}(N)^i = 
S(i)\otimes_kN_i$. For every $j$ one has: 
\[
\text{H}^i(\text{F}(N))_j\simeq \text{Hom}_{\text{K}^b(\mathbb P)}
(\text{T}^{-i}{\mathcal O}_{\mathbb P},\text{L}(N)(j))
\]
hence it remains to show that for $j\geq -i$:
\[
\text{Hom}_{\text{K}^b(\mathbb P)}(\text{T}^{-i}{\mathcal O}_{\mathbb P},
\text{L}(N)(j))\overset{\sim}{\longrightarrow} 
\text{Hom}_{\text{D}^b(\mathbb P)}(\text{T}^{-i}{\mathcal O}_{\mathbb P},
\text{L}(N)(j))
\]
or, equivalently, that:
\[
\text{Hom}_{\text{K}^b(\mathbb P)}(\text{L}(\underline{k}(-i)),
\text{T}^{i+j}\text{L}(N(-i-j)))\overset{\sim}{\longrightarrow} 
\text{Hom}_{\text{D}^b(\mathbb P)}(\text{L}(\underline{k}(-i)),
\text{T}^{i+j}\text{L}(N(-i-j))).
\]
For $j > -i$ this follows from Lemma 4.2.. For $j = -i$, the above morphism 
can be identified with the morphism:
\[
\text{Hom}_{\Lambda\text{-mod}}(\underline{k}(-i),N)\longrightarrow 
\text{Hom}_{\text{D}^b(\mathbb P)}(\text{L}(\underline{k}(-i)),
\text{L}(N)).
\]
By Theorem 4.4., the last morphism is surjective and its kernel consists of 
the composite morphisms $\underline{k}(-i)\rightarrow P\rightarrow N$ with 
$P$ a free object of $\Lambda$-mod. But the image of $\underline{k}(-i)
\rightarrow P$ must lie in $\text{soc}(P) = \text{soc}(\Lambda)\cdot P$. 
Since $N$ is annihilated by $\text{soc}(\Lambda)$, any such composite 
morphism must be 0.
\end{proof}
\vskip3mm

{\bf 4.7. Definition.}\quad 
Let $N$ be an object of $\Lambda$-mod annihilated by $\text{soc}(\Lambda)$. 
Let $P^{\bullet}$ (resp., $P^{\prime \bullet}$) be a minimal free resolution 
of $N$ (resp., $N^{\ast}$) in $\Lambda$-mod. By concatenating the complexes 
$P^{\prime \bullet}$ and $\text{T}^{-1}(P^{\bullet \ast})$ using the composite 
morphism $P^{\prime 0}\twoheadrightarrow N^{\ast}\hookrightarrow 
P^{0\ast}$ one gets an {\it acyclic complex} which is {\it minimal} 
(since $N^{\ast}$ is annihilated by $\text{soc}(\Lambda)$, 
the image of the above composite 
morphism is contained in ${\Lambda}_+\cdot P^{0\ast}$). The $k$-vector 
space dual $I^{\bullet}$ of this complex (see Definition 2.5.) is called 
a {\it Tate resolution} of $N$.
\vskip3mm 

The next theorem, which is one of the main results of the paper of 
Eisenbud, Fl\o ystad and Schreyer \cite{efs}, is a direct consequence of 
Corollary 4.6..

\vskip3mm
{\bf 4.8. Theorem.} (\cite{efs}, Theorem 4.1.)\quad 
\textit{If} ${\mathcal F}^{\bullet}$ \textit{and} $N$ \textit{are as in 
Corollary} 4.5. \textit{and if} $I^{\bullet}$ \textit{is a Tate resolution of} 
$N$ \textit{then the linear part of} $I^{\bullet}$ \textit{is isomorphic to} 
$\bigoplus_{i\in {\mathbb Z}}\text{T}^{-i}\text{G}({\mathbb H}^i_{\ast}
{\mathcal F}^{\bullet})$. 
\vskip3mm

\begin{proof}
Let $Z^{-m} := \text{Ker}(I^{-m}\rightarrow I^{-m+1})$, $m > 0$. 
As we shown in the proof of Lemma 3.8., the linear part of 
${\sigma}^{\geq -m}I^{\bullet}$ is isomorphic to  
$\bigoplus_{i\in {\mathbb Z}}\text{T}^{-i}\text{G}(
\text{H}^{i+m}(\text{F}(Z^{-m})))$.  

Now, by definition, $N \simeq \text{Coker}(I^{-2}\rightarrow I^{-1})$ and, 
since $I^{\bullet}$ is acyclic, $N \simeq \text{Ker}(I^0\rightarrow I^1)$. 
It follows, from 
observation (2) in the second part of the proof of Theorem 4.4., that  
$\text{L}(Z^{-m}) \simeq \text{T}^{-m}\text{L}(N) \simeq \text{T}^{-m} 
{\mathcal F}^{\bullet}$ in $\text{D}^b(\mathbb P)$. Moreover, since $Z^{-m}$ 
is contained in ${\Lambda}_+\cdot I^{-m}$, it is annihilated by 
$\text{soc}(\Lambda)$. Corollary 4.6. implies, now, that 
$\text{H}^{i+m}(\text{F}(Z^{-m}))\simeq 
\bigoplus_{j\geq -i-m}{\mathbb H}^{i+m}((\text{T}^{-m}{\mathcal F}^{\bullet})
(j)) = \bigoplus_{j\geq -i-m}{\mathbb H}^i({\mathcal F}^{\bullet}(j))$. 
Taking into account what have been recalled in the first paragraph, 
one deduces 
that the linear part of ${\sigma}^{\geq -m}I^{\bullet}$ is isomorphic to 
$\bigoplus_{i\in {\mathbb Z}}\text{T}^{-i}\text{G}
(\bigoplus_{j\geq -i-m}{\mathbb H}^i({\mathcal F}^{\bullet}(j)))$. 
Finally, letting $m\rightarrow \infty$ one gets the desired conclusion.
\end{proof}
\vskip3mm

{\bf 4.9. Theorem.}\quad 
\textit{Let} $\mathcal F$ \textit{be a coherent sheaf on} ${\mathbb P}^n$ 
\textit{with} $\text{H}^0{\mathcal F}(-t) = 0$ \textit{for} $t >> 0$, 
\textit{let} $M := \text{H}^0_{\ast}{\mathcal F}$ \textit{and let} 
$0 \rightarrow L^{-n+1} \rightarrow \cdots \rightarrow L^0 \rightarrow M 
\rightarrow 0$ \textit{be a minimal free resolution of} $M$ \textit{in} 
$S$-mod. \textit{Let} $N\in \text{Ob}(\Lambda\text{-mod})$ \textit{be as in 
Corollary} 4.5. \textit{and let} $I^{\bullet}$ \textit{be a Tate resolution 
of} $N$. \textit{Then}:\\
\hspace*{3mm} (a) $I^{\bullet}/F_0I^{\bullet}$ \textit{is a contraction of} 
$\text{T}^{-n}\text{G}((L^{\bullet \vee}\otimes_S{\omega}_S)^{\ast})$.\\
\hspace*{3mm} (b) \textit{The Horrocks-Trautmann complex corresponding to} 
$\mathcal F$ \textit{via the equivalence of categories from Theorem} 3.9. 
\textit{is isomorphic to} $F_{n-1}I^{\bullet}/F_0I^{\bullet}$.
\vskip3mm

\begin{proof}
(a) Choose $m\in {\mathbb Z}$ such that $\text{H}^i({\mathcal F}(j)) = 0$, 
$\forall i > 0$, $\forall j\geq m-i$. Since, as a 
consequence of Theorem 4.8., $I^p \simeq \bigoplus_{i=0}^n\text{H}^i
({\mathcal F}(p-i))\otimes_k\bigwedge (V^{\ast})(p-i)$, $\forall p\in 
{\mathbb Z}$, one sees that $I^p = F_0I^p$ for $p\geq m$, hence $I^{\bullet}/
F_0I^{\bullet} = ({\sigma}^{<m}I^{\bullet})/F_0({\sigma}^{<m}I^{\bullet})$. 

Now, let $Z^m := \text{Ker}(I^m\rightarrow I^{m+1})$. 
$\text{T}^{-1}({\sigma}^{<m}I^{\bullet})$ is a minimal (left) free resolution 
of $\text{T}^{-m}Z^m$. One deduces, from Corollary 2.12.,  
that $\text{T}^{-1}({\sigma}^{<m}I^{\bullet})$ is a contraction of 
$\text{T}^{-n-1}\text{G}((\text{F}(\text{T}^{-m}Z^m)^{\vee}
\otimes_S{\omega}_S)^{\ast})$, hence ${\sigma}^{<m}I^{\bullet}$ is a 
contraction of $\text{T}^{-n}\text{G}((\text{F}(\text{T}^{-m}Z^m)^{\vee}
\otimes_S{\omega}_S)^{\ast})$. 

Let $F^{\bullet} := \text{F}(\text{T}^{-m}Z^m) = \text{T}^{-m}\text{F}(Z^m)$. 
By observation (2) in the last part of the proof of Theorem 4.4., 
$\text{L}(Z^m)\simeq \text{T}^m\text{L}(N)\simeq \text{T}^m{\mathcal F}$ in 
$\text{D}^b(\mathbb P)$. By Corollary 4.6., $\text{H}^i(F^{\bullet}) = 
\text{H}^{i-m}(\text{F}(Z^m)) \simeq \bigoplus_{j\geq m-i}{\mathbb H}^{i-m}
((\text{T}^m{\mathcal F})(j)) = \bigoplus_{j\geq m-i}\text{H}^i
({\mathcal F}(j))$. One deduces that $F^{\bullet}$ is a 
minimal free resolution of $M^{\prime} := \bigoplus_{j\geq m}\text{H}^0
({\mathcal F}(j))$ in $S$-mod. 

The inclusion 
$M^{\prime} \hookrightarrow M$ lifts to a morphism of complexes $F^{\bullet} 
\rightarrow L^{\bullet}$. Since $M/M^{\prime}$ is of finite length, 
$\text{Ext}^i_S(M/M^{\prime},{\omega}_S) = 0$ for $i \neq n+1$, hence 
$\text{H}^i(L^{\bullet \vee}\otimes_S{\omega}_S) \overset{\sim}{\rightarrow} 
\text{H}^i(F^{\bullet \vee}\otimes_S{\omega}_S)$ for $i < n$. One deduces a 
quasi-isomorphism 
\[
L^{\bullet \vee}\otimes_S{\omega}_S \rightarrow 
{\tau}^{<n}(F^{\bullet \vee}\otimes_S{\omega}_S), 
\]
whence a quasi-isomorphism 
${\tau}^{>-n}((F^{\bullet \vee}\otimes_S{\omega}_S)^{\ast})\rightarrow 
(L^{\bullet \vee}\otimes_S{\omega}_S)^{\ast}$ and then a quasi-isomorphism 
\[
{\tau}^{>0}(\text{T}^{-n}((F^{\bullet \vee}\otimes_S{\omega}_S)^{\ast})) 
\rightarrow \text{T}^{-n}((L^{\bullet \vee}\otimes_S{\omega}_S)^{\ast}). 
\]
It follows now, from the last part of Lemma 2.9.(b), that 
$({\sigma}^{<m}I^{\bullet})/F_0({\sigma}^{<m}I^{\bullet})$ is a contraction 
of $\text{T}^{-n}\text{G}((L^{\bullet \vee}\otimes_S{\omega}_S)^{\ast})$. 

(b) Let $0\rightarrow L^{\prime -n+1}\rightarrow \cdots \rightarrow 
L^{\prime 0}\rightarrow M^{\vee}\rightarrow 0$ be a minimal free resolution 
of $M^{\vee}$ in $S$-mod. As in the last part of Definition 3.2., one can 
construct from $L^{\bullet}$ and $L^{\prime \bullet}$ a Horrocks resolution 
$K^{\bullet}$ of $M$. One has an exact sequence:
\[
0\rightarrow \text{T}^{-1}(L^{\bullet \vee})\rightarrow K^{\bullet \vee}
\rightarrow L^{\prime \bullet}\rightarrow 0.
\tag{*}
\]
The Horrocks-Trautmann complex corresponding to $\mathcal F$ is obtained as 
follows: one considers a complex $N^{\bullet}\in \text{Ob}\, \text{C}^b
(\Lambda \text{-mod})$ such that $\text{F}(N^{\bullet}) \simeq K^{\bullet}$ in 
$\text{K}^b(\mathcal P)$ and then one takes a minimal (left) free resolution 
$G^{\bullet}$ of $N^{\bullet}$ in $\Lambda$-mod. By Corollary 2.12., 
$G^{\bullet}$ is a contraction of $\text{T}^{-n-1}\text{G}((\text{F}
(N^{\bullet})^{\vee}\otimes_S{\omega}_S)^{\ast})$ hence a contraction of 
$\text{T}^{-n-1}\text{G}((K^{\bullet \vee}\otimes_S{\omega}_S)^{\ast})$. 

Using the exact sequence (*), one gets a quasi-isomorphism: 
\[
{\tau}^{\geq 2}(\text{T}^{-1}(L^{\bullet \vee}\otimes_S{\omega}_S)) 
\longrightarrow K^{\bullet \vee}\otimes_S{\omega}_S,
\]
whence a quasi-isomorphism $(K^{\bullet \vee}\otimes_S{\omega}_S)^{\ast} 
\rightarrow {\tau}^{\leq -2}(\text{T}((L^{\bullet \vee}\otimes_S{\omega}_S
)^{\ast}))$ and then a quasi-isomorphism: 
\[
\text{T}^{-n-1}((K^{\bullet \vee}\otimes_S{\omega}_S)^{\ast})\longrightarrow 
{\tau}^{\leq n-1}(\text{T}^{-n}((L^{\bullet \vee}\otimes_S{\omega}_S)^{\ast}
)).
\]
One deduces, now, from (a) and from the last part of Lemma 2.9.(b), that the 
contraction $G^{\bullet}$ of $\text{T}^{-n-1}\text{G}((K^{\bullet \vee}
\otimes_S{\omega}_S)^{\ast})$ is isomorphic to $F_{n-1}I^{\bullet}/
F_0I^{\bullet}$.     
\end{proof}    

\section*{Appendix A : Cancellation of terms in a complex}

We work in an abelian category $\mathcal A$. 

\vskip3mm
{\bf A.1. Definition.}\quad
If $X^{\bullet}$ and $Y^{\bullet}$ are complexes in  
$\mathcal A$ then, 
according to Eilenberg and MacLane \cite{emc}, a {\it contraction} of 
$X^{\bullet}$ onto $Y^{\bullet}$ is a triple $(f,g,h)$, where 
$f:X^{\bullet}\rightarrow Y^{\bullet}$, $g:Y^{\bullet}\rightarrow 
X^{\bullet}$ are morphisms of complexes and $h\in {\text{Hom}}^{-1}
(X^{\bullet},X^{\bullet})$ is a homotopy operator satisfying: 
\[
(\text{i})\  fg={\text{id}}_Y,\  \  (\text{ii})\  {\text{id}}_X-gf=d_Xh+hd_X
\]
and the {\it side conditions}: 
\[
(\text{iii})\  fh=0,\  \  (\text{iv})\  hg=0,\  \  (\text{v})\  h^2=0.
\]
The side conditions do not restrict generality. Indeed, as remarked by 
Lambe and Stasheff \cite{lst}, (2.1.), if $(f,g,h)$ satisfies (i)-(ii), 
if one puts $\phi := \text{id}_X - gf$ and $h^{\prime}:=\phi h\phi $, then 
$(f,g,h^{\prime})$ satisfies (i)-(iv) ($d_X\phi =\phi d_X$, $d_Xh + hd_X = 
\phi $, and ${\phi}^2 = \phi $). 
Moreover, if  $h^{\prime \prime}:=h^{\prime}d_Xh^{\prime}$ then 
$(f,g,h^{\prime \prime})$ satisfies (i)-(v) ($h^{\prime} = h^{\prime}\phi = 
\phi h^{\prime}$, and $\phi = d_Xh^{\prime} + h^{\prime}d_X$). 

Notice that (i) implies that $X^{\bullet} = \text{Im}g\oplus 
\text{Ker}f$, (iv) implies that $h$ vanishes on $\text{Im}g$, 
(iii) implies that $h$ maps $X^{\bullet}$ into $\text{Ker}f$, and 
(ii) implies that $h\, \vert \, \text{Ker}f \rightarrow \text{Ker}f$ 
realizes a homotopy $\text{id}_{\text{Ker}f} \sim 0$.   
\vskip3mm

{\bf A.2. Example.}\quad
Let $X^{\bullet}$ be a complex in $\mathcal A$ and assume that, for every 
$p\in {\mathbb Z}$, one has a decomposition $X^p=V^p\oplus W^p\oplus Y^p$ 
such that the component $d^p_{wv}:V^p\rightarrow W^{p+1}$ of the 
differential $d^p_X:X^p\rightarrow X^{p+1}$ is an {\it isomorphism}. 
Consider, for every $p\in {\mathbb Z}$, the following morphisms:
\begin{gather*}
d^p_Y=d^p_{yy}-d^p_{yv}(d^p_{wv})^{-1}d^p_{wy}: Y^p\longrightarrow Y^{p+1},\\
f^p=(0\  ,\  -d^{p-1}_{yv}(d^{p-1}_{wv})^{-1}\  ,\  {\text{id}}_{Y^p}): X^p
\longrightarrow Y^p,\\
g^p=(-(d^p_{wv})^{-1}d^p_{wy}\  ,\  0\  ,\  {\text{id}}_{Y^p})^{\text{trsp}}: 
Y^p\longrightarrow X^p,
\end{gather*}
and $h^p: X^p\rightarrow X^{p-1}$ defined by the $3\times 3$ matrix whose 
unique non-zero entry is $(d^{p-1}_{wv})^{-1}:W^p\rightarrow V^{p-1}$. 
Then $Y^{\bullet}:=(Y^p,d^p_Y)_{p\in {\mathbb Z}}$ is a complex, 
$f:=(f^p)_{p\in {\mathbb Z}}$ and $g:=(g^p)_{p\in {\mathbb Z}}$ are 
morphisms of complexes and $(f,g,h:=(h^p)_{p\in {\mathbb Z}})$ 
is a contraction of $X^{\bullet}$ onto $Y^{\bullet}$. 
\vskip3mm

{\bf A.3. Example.}\quad
An important particular case of the previous example is that of a 
{\it splitting complex} $X^{\bullet}$. This means that, $\forall p\in 
{\mathbb Z}$, $Z^p:=\text{Ker}d^p_X$ and $B^p:=\text{Im}d^{p-1}_X$ are 
direct summands of $X^p$. Choose decompositions : $X^p=V^p\oplus Z^p$ and 
$Z^p=B^p\oplus H^p$, hence $X^p=V^p\oplus B^p\oplus H^p$. The differential 
$d^p_X$ vanishes on $Z^p$ and maps $V^p$ isomorphically onto $B^{p+1}$. Let 
$h^p:X^p\rightarrow X^{p-1}$ be the morphism defined by the $3\times 3$ 
matrix whose unique non-zero term is the inverse of 
$V^{p-1}\overset{\sim}{\rightarrow} B^p$. Consider the complex  
$H^{\bullet}:=(H^p,0)_{p\in {\mathbb Z}}$. Then the projection (corresponding  
to the above decompositions) $\pi :X^{\bullet}
\twoheadrightarrow H^{\bullet}$ and the inclusion $u:H^{\bullet}
\hookrightarrow X^{\bullet}$ are morphisms of complexes and $(\pi ,u,h)$ is 
a contraction of $X^{\bullet}$ onto $H^{\bullet}$.

Remark that the composite morphisms $H^p\hookrightarrow Z^p
\twoheadrightarrow {\text{H}}^p(X^{\bullet})$, $p\in {\mathbb Z}$, define an 
isomorphism of complexes $H^{\bullet}\overset{\sim}{\rightarrow} 
{\text{H}}^{\bullet}(X^{\bullet}):=
({\text{H}}^p(X^{\bullet}),0)_{p\in {\mathbb Z}}$. 

Remark, also, that, \textit{conversely}, \textit{if a complex} 
$X^{\bullet}$ \textit{is homotopically equivalent to a complex} 
$H^{\bullet}$ \textit{with} $d_H = 0$ \textit{then} $X^{\bullet}$ 
\textit{is a splitting complex}. 

Indeed, consider morphisms of complexes $f : X^{\bullet} \rightarrow 
H^{\bullet}$ and $g : H^{\bullet} \rightarrow X^{\bullet}$ such that 
$gf \sim \text{id}_X$ and $fg \sim \text{id}_H$. Choose a homotopy operator 
$h \in \text{Hom}^{-1}(H^{\bullet},H^{\bullet})$ such that $\text{id}_H - fg = 
d_Hh + hd_H$. Since $d_H = 0$, it follows that $fg = \text{id}_H$, hence 
$X^{\bullet} \simeq X^{\prime \bullet} \oplus H^{\bullet}$, where 
$X^{\prime \bullet} = \text{Ker}f$. Choose, now, a homotopy operator 
$k \in \text{Hom}^{-1}(X^{\bullet},X^{\bullet})$ such that 
$\text{id}_X - gf = d_Xk + kd_X$ and $fk = 0$ 
(see the argument of Lambe and Stasheff reproduced in (A.1.)). Then 
$k$ maps $X^{\bullet}$ into $\text{Ker}f = X^{\prime \bullet}$ and 
induces a homotopy operator $h^{\prime}$ on $X^{\prime \bullet}$ such that 
$\text{id}_{X^{\prime}} = d_{X^{\prime}}h^{\prime} + h^{\prime}d_{X^{\prime}}$. 
In particular, $X^{\prime \bullet}$ is acyclic, hence $B^{\prime p} := 
\text{Im}d^{p-1}_{X^{\prime}} = \text{Ker}d^p_{X^{\prime}} =: Z^{\prime p}$. 
Now, $d^{p-1}_{X^{\prime}} \circ h^{\prime p}\, \vert \, X^{\prime p} 
\rightarrow B^{\prime p}$ is a left inverse for the inclusion 
$Z^{\prime p} \hookrightarrow X^{\prime p}$.   
\vskip3mm

The next result is known in the literature as the ``Basic Perturbation  
Lemma''. In its more practical form (A.6.) below, it appears implicitly in 
Shih \cite{shi} and explicitly in R. Brown \cite{bro} and Gugenheim 
\cite{gug}. Its more general variant (A.4.) was proved by Barnes and Lambe 
\cite{bar}. We include here a different proof of this variant. 

\vskip3mm
{\bf A.4. Basic Perturbation Lemma.}\quad
\textit{Let} $(f,g,h)$ 
\textit{be a contraction of a complex} $X^{\bullet}$ 
\textit{onto a complex} $Y^{\bullet}$. 
\textit{Let} ${\widehat d}_X=d_X+d^{\prime}_X\in 
\text{Hom}^1(X^{\bullet},X^{\bullet})$ \textit{be a} ``\textit{perturbation}'' 
\textit{of} $d_X$ 
(\textit{which means that} ${\widehat d}_X\circ {\widehat d}_X=0$, 
\textit{i.e.}, \textit{that} ${\widehat X}^{\bullet}:= 
(X^p,{\widehat d}^{\  p}_X)_{p\in {\mathbb Z}}$ \textit{is a complex}). 
\textit{If}  
$\text{id}_X+hd^{\prime}_X$ \textit{is an invertible element of the ring}  
$\text{Hom}^0(X^{\bullet},X^{\bullet})$ 
\textit{then there exist a perturbation}  
${\widehat d}_Y=d_Y+d^{\prime}_Y$ \textit{of} $d_Y$ 
\textit{and a contraction}  
$({\widehat f},{\widehat g},{\widehat h})$ \textit{of the complex}  
${\widehat X}^{\bullet}$ \textit{onto the complex} 
${\widehat Y}^{\bullet}:=
(Y^p,{\widehat d}^{\  p}_Y)_{p\in {\mathbb Z}}$. 
\vskip3mm
 
\begin{proof}
We consider, firstly, the {\it particular case} where $d^{\prime}_X=g\phi $ 
for some $\phi \in {\text{Hom}}^1(X^{\bullet},Y^{\bullet})$. In this case:
\[
{\widehat d}_Xg=d_Xg+g\phi g=gd_Y+g\phi g=g(d_Y+\phi g),
\]
hence, putting ${\widehat d}_Y:=d_Y+\phi g$, one gets that 
${\widehat d}_Xg=g{\widehat d}_Y$. It follows that 
\[
{\widehat d}_Y\circ 
{\widehat d}_Y=f\circ g\circ {\widehat d}_Y\circ {\widehat d}_Y=
f\circ {\widehat d}_X\circ {\widehat d}_X\circ g=0,
\]
i.e., ${\widehat d}_Y$ is a perturbation of $d_Y$. 

We look, now, for a perturbation ${\widehat f}=f+f^{\prime}$ of $f$ such 
that:
\begin{gather}
{\widehat f}g={\text{id}}_Y,\  \  {\text{id}}_X-g{\widehat f}=
{\widehat d}_Xh+h{\widehat d}_X\tag{*}.
\end{gather}
This system of equations is equivalent to:
\begin{gather}
f^{\prime}g=0,\  \  -gf^{\prime}=g\phi h\tag{**}
\end{gather}
(because $hg=0$) hence it has the solution $f^{\prime}:=-\phi h$. Now, by 
(*), ${\widehat f}g={\text{id}}_Y$ and $g{\widehat f}$ is an endomorphism 
of the complex ${\widehat X}^{\bullet}$ hence:
\[
{\widehat d}_Y{\widehat f}-{\widehat f}{\widehat d}_X=
{\widehat f}g({\widehat d}_Y{\widehat f}-{\widehat f}{\widehat d}_X)=
{\widehat f}({\widehat d}_Xg{\widehat f}-{\widehat d}_Xg{\widehat f})=0,
\]
i.e., ${\widehat f}$ is a morphism of complexes from ${\widehat X}^{\bullet}$ 
to ${\widehat Y}^{\bullet}$. Moreover, ${\widehat f}h=-\phi h^2=0$. 
Consequently, $({\widehat f},g,h)$ is a contraction of 
${\widehat X}^{\bullet}$ onto ${\widehat Y}^{\bullet}$. 

The {\it general case} can be reduced to the particular case we have just 
treated as follows : $\alpha :={\text{id}}_X+hd^{\prime}_X\in 
{\text{Hom}}^0(X^{\bullet},X^{\bullet})$ maps isomorphically the complex 
${\widehat X}^{\bullet}$ onto the complex ${\widetilde X}^{\bullet}=
(X^p,{\widetilde d}^{\  p}_X)_{p\in {\mathbb Z}}$, where ${\widetilde d}_X:=
\alpha {\widehat d}_X{\alpha}^{-1}$. Using the fact that 
$d_Xd_X^{\prime} + d_X^{\prime}{\widehat d}_X = 
{\widehat d}_X\circ {\widehat d}_X = 0$ 
and the relation (ii) from (A.1.) one checks easily that:
\[
\alpha {\widehat d}_X=d_X\alpha +gfd^{\prime}_X
\]
hence ${\widetilde d}_X=\alpha {\widehat d}_X{\alpha}^{-1}=d_X+
gfd^{\prime}_X{\alpha}^{-1}$ is a perturbation of $d_X$ with perturbation 
term $g\phi $, where $\phi =fd^{\prime}_X{\alpha}^{-1}$. It follows, from 
the particular case, that there exist perturbations:
\[
{\widehat d}_Y=d_Y+fd^{\prime}_X{\alpha}^{-1}g,\  \  
{\widetilde f}=f-fd^{\prime}_X{\alpha}^{-1}h
\]
such that $({\widetilde f},g,h)$ is a contraction of the complex 
${\widetilde X}^{\bullet}$ onto the complex ${\widehat Y}^{\bullet}$. 

One can now take: ${\widehat f}={\widetilde f}\alpha ={\widetilde f}$ 
(because $fh=0$ and $h^2=0$), ${\widehat g}={\alpha}^{-1}g$ and 
${\widehat h}={\alpha}^{-1}h\alpha ={\alpha}^{-1}h$ (because $h^2=0$).  
\end{proof} 

\vskip3mm
{\bf A.5. Remark.}\quad 
\textit{Under the hypothesis of} (A.4.), \textit{let} $U^{\bullet} := 
\text{Ker}f$. \textit{The sequence}  
$0\rightarrow U^{\bullet}\overset{u}{\longrightarrow} X^{\bullet} 
\overset{f}{\longrightarrow} Y^{\bullet}\rightarrow 0$ 
\textit{is split exact} (\textit{because} $fg = \text{id}_Y$), 
\textit{hence} $X^{\bullet} \simeq Y^{\bullet} \oplus U^{\bullet}$. 
\textit{Then}, \textit{for the complex} ${\widehat Y}^{\bullet}$ 
\textit{obtained in the proof of} (A.4.), 
\textit{one has}: ${\widehat X}^{\bullet} \simeq {\widehat Y}^{\bullet} 
\oplus U^{\bullet}$. 
\vskip3mm

\begin{proof} 
Since $f(\text{id}_X - gf) = 0$, it follows that there exists a morphism 
of complexes $v : X^{\bullet} \rightarrow U^{\bullet}$ such that 
$\text{id}_X - gf = uv$. One deduces that $vu = \text{id}_U$ and that 
$vg = 0$, hence the sequence of complexes  
$0\rightarrow Y^{\bullet}\overset{g}{\longrightarrow} X^{\bullet} 
\overset{v}{\longrightarrow} U^{\bullet}\rightarrow 0$ 
is (split) exact. 

Now, using the notation from the last part of the proof of (A.4.), one 
has: 
\[
v{\widetilde d}_X = vd_X + vg\phi = vd_X = d_Uv
\]
hence $v$ is a morphism of complexes: ${\widetilde X}^{\bullet} 
\rightarrow U^{\bullet}$. The short exact sequence of complexes: 
\[
0\rightarrow {\widehat Y}^{\bullet}\overset{g}{\longrightarrow} 
{\widetilde X}^{\bullet}\overset{v}{\longrightarrow} U^{\bullet} 
\rightarrow 0
\]
is split exact (because ${\widetilde f}g = \text{id}_{\widehat Y}$), hence 
${\widetilde X}^{\bullet} \simeq {\widehat Y}^{\bullet} \oplus U^{\bullet}$. 
But one has an isomorphism of complexes $\alpha : {\widehat X}^{\bullet} 
\overset{\sim}{\rightarrow} {\widetilde X}^{\bullet}$. 
\end{proof} 

\vskip3mm
{\bf A.6. The classical variant.}\quad
In practice, one checks the fact that ${\text{id}}_X+hd^{\prime}_X$ is 
invertible by verifying that $hd^{\prime}_X$ is {\it locally nilpotent}, i.e., 
that ${\bigcup}_{i\geq 1}\text{Ker}(hd^{\prime}_X)^i=X^{\bullet}$. In this 
case, the inverse of ${\text{id}}_X+hd^{\prime}_X$ is 
${\text{id}}_X+\sum_{i\geq 1}(-1)^i(hd^{\prime}_X)^i$ and the proof of 
(A.4.) gives the following explicit formulae:
\begin{gather*}
{\widehat d}_Y=d_Y+fd^{\prime}_Xg+
\sum_{i\geq 1}(-1)^ifd^{\prime}_X(hd^{\prime}_X)^ig,
\  \  {\widehat f}=f+\sum_{i\geq 1}(-1)^if(d^{\prime}_Xh)^i,\\ 
{\widehat g}=g+\sum_{i\geq 1}(-1)^i(hd^{\prime}_X)^ig,\  \  
{\widehat h}=h+\sum_{i\geq 1}(-1)^i(hd^{\prime}_X)^ih.
\end{gather*}
\vskip3mm

{\bf A.7. The case of a double complex.}\quad 
Let $X^{\bullet \bullet}$ be a double complex with (commuting) differentials 
$d^{\prime}_X$ and $d^{\prime \prime}_X$. We denote by 
${\text{H}}_{\text{I}}(X^{\bullet \bullet})$ 
the double complex with terms 
${\text{H}}_{\text{I}}^{pq}(X^{\bullet \bullet})
:=\text{Ker}d^{\prime pq}_X/
\text{Im}d^{\prime p-1,q}_X$, with 
$d^{\prime}_{{\text{H}}_{\text{I}}}=0$ and with 
$d^{\prime \prime}_{{\text{H}}_{\text{I}}}$ 
induced by $d^{\prime \prime}_X$. We also recall 
the following notation: for $m\in {\mathbb Z}$, 
${\tau}_{\text{I}}^{\leq m}X^{\bullet \bullet}$ is the double subcomplex of 
$X^{\bullet \bullet}$ whose $(p,q)$ term is $X^{pq}$ for $p<m$, 
$\text{Ker}d^{\prime mq}_X$ for $p=m$, and 0 for $p>m$. One defines, 
similarly, a quotient double complex ${\tau}^{>m}_{\text{I}}X^{\bullet \bullet}$ 
of $X^{\bullet \bullet}$.  

The following result, which is a particular case of (A.6.), is stated and 
proved in Eisenbud et al. \cite{efs}, (3.5), and it is a key technical 
point of that paper. 

\vskip3mm
{\bf Lemma.}\quad 
\textit{Assume that the double complex} 
$X^{\bullet \bullet}$ \textit{satisfies the following finiteness condition}:  
$\forall m\in {\mathbb Z}$, $X^{p,m-p}=0$ \textit{for}  
$p<<0$. \textit{If all the rows} 
$X^{\bullet ,q}:=(X^{pq},d^{\prime pq}_X)_{p\in {\mathbb Z}}$, 
$q\in {\mathbb Z}$, 
\textit{of} $X^{\bullet \bullet}$ \textit{split}  
(\textit{see} (A.3.)) \textit{then there exists a contraction of} 
$\text{tot}(X^{\bullet \bullet})$ \textit{onto a complex} 
$Y^{\bullet}$, \textit{endowed with an increasing filtration} 
$(F_mY^{\bullet})_{m\in {\mathbb Z}}$ \textit{by subcomplexes, such that}:
\begin{gather}
Y^n=\bigoplus_{p+q=n}\text{H}^{pq}_{\text{I}}(X^{\bullet \bullet}),\   
\forall n\in {\mathbb Z},\tag{1}\\  
F_mY^n=\bigoplus_{\substack{p+q=n\\p\leq m}}\text{H}^{pq}_{\text{I}}
(X^{\bullet \bullet}),\  \forall m,n\in {\mathbb Z},\tag{2}\\   
\text{gr}_F(Y^{\bullet})=\text{tot}(\text{H}_{\text{I}}
(X^{\bullet \bullet})).
\tag{3}
\end{gather}
\textit{Moreover, this contraction can be chosen in such a way that, for all}  
$m\in {\mathbb Z}$, \textit{it induces a contraction of} $\text{tot}
({\tau}_{\text{I}}^{\leq m}X^{\bullet \bullet})$ \textit{onto}  
$F_mY^{\bullet}$ \textit{and of}  
$\text{tot}({\tau}^{>m}_IX^{\bullet \bullet})$ \textit{onto} 
$Y^{\bullet}/F_mY^{\bullet}$.  
\vskip3mm

\begin{proof}
Recall that the differential of $\text{tot}(X^{\bullet \bullet})$ is 
$d^{\prime}_X+{\delta}^{\prime \prime}_X$, where ${\delta}^{\prime \prime}_X 
\vert X^{pq}:=(-1)^pd^{\prime \prime pq}_X$. 
Let $X_{\text{I}}^{\bullet \bullet}$ 
be the double complex with the same terms as $X^{\bullet \bullet}$,  
with $d_{X_{\text{I}}}^{\prime}=d_X^{\prime}$, and with 
$d_{X_{\text{I}}}^{\prime \prime}=0$. (A.3.) provides a contraction 
$(\pi ,u,h)$ of $\text{tot}(X_{\text{I}}^{\bullet \bullet})$ onto a complex 
with terms $Y^n$ given by the formula (1) from the statement and with the 
differential equal to 0. One may assume that the homotopy operator $h$ maps 
$X^{pq}$ into $X^{p-1,q}$, $\forall p,q\in {\mathbb Z}$. The differential of 
$\text{tot}(X^{\bullet \bullet})$ is a perturbation of $d^{\prime}_X$ and 
the finiteness condition from the statement implies that 
$h{\delta}^{\prime \prime}_X$ is 
locally nilpotent. (A.6.) produces now a contraction $({\widehat \pi},
{\widehat u},{\widehat h})$ of $\text{tot}(X^{\bullet \bullet})$ onto a 
complex $Y^{\bullet}$ with terms given by formula (1) from the statement 
and with differential:
\[
d_Y=\pi {\delta}^{\prime \prime}_Xu+\sum_{i\geq 1}(-1)^i\pi 
{\delta}^{\prime \prime}_X(h{\delta}^{\prime \prime}_X)^iu. 
\]
The explicit formulae from (A.6.) allows one now to check easily the other 
assertions from the lemma.
\end{proof}

\section*{Appendix B : A comparison lemma}

\vskip3mm
{\bf B.1. Lemma.}\quad 
\textit{Let} $\mathcal{C}$, $\mathcal{D}$ \textit{be triangulated categories 
and} $\Phi : \mathcal{C}\rightarrow \mathcal{D}$ \textit{an additive functor 
commuting with the translation functors and sending distinguished triangles 
to distinguished triangles}. \textit{Let} $X$, $Y$ \textit{be two objects of} 
$\mathcal{C}$ \textit{endowed with} ``\textit{decreasing filtrations}'', 
\textit{i.e., with sequences of morphisms}: 
\[
\cdots \rightarrow F^{i+1}X\rightarrow F^iX\rightarrow \cdots ,\  \  
\cdots \rightarrow F^{i+1}Y\rightarrow F^iY\rightarrow \cdots 
\]
\textit{such that} $F^iX = X$, $F^iY = Y$ \textit{for} $i << 0$ \textit{and} 
$F^iX = 0$, $F^iY = 0$ \textit{for} $i >> 0$, \textit{and with the} 
``\textit{successive quotients}'' \textit{replaced by distinguished 
triangles}: 
\[
F^{i+1}X\rightarrow F^iX\rightarrow X^i\rightarrow \text{T}F^{i+1}X,\  \  
F^{i+1}Y\rightarrow F^iY\rightarrow Y^i\rightarrow \text{T}F^{i+1}Y.
\]
\hspace*{3mm} (a) \textit{If} $\text{Hom}_{\mathcal{C}}(X^i,Y^j)
\rightarrow \text{Hom}_{\mathcal{D}}(\Phi (X^i),\Phi (Y^j))$ \textit{is 
surjective and} $\text{Hom}_{\mathcal{C}}(X^i,\text{T}Y^j)\rightarrow 
\text{Hom}_{\mathcal{D}}(\Phi (X^i),\Phi (\text{T}Y^j))$ 
\textit{is injective}, $\forall i,j$, \textit{then} 
$\text{Hom}_{\mathcal{C}}(X,Y)\rightarrow \text{Hom}_{\mathcal{D}}
(\Phi (X),\Phi (Y))$ \textit{is surjective}.\\ 
\hspace*{3mm} (b) \textit{If} $\text{Hom}_{\mathcal{C}}(X^i,Y^j)\rightarrow 
\text{Hom}_{\mathcal{D}}(\Phi (X^i),\Phi (Y^j))$ \textit{is injective and} 
$\text{Hom}_{\mathcal{C}}(\text{T}X^i,Y^j)\rightarrow 
\text{Hom}_{\mathcal{D}}(\Phi (\text{T}X^i),\Phi (Y^j))$ \textit{is 
surjective}, $\forall i,j$, \textit{then} $\text{Hom}_{\mathcal{C}}(X,Y)
\rightarrow \text{Hom}_{\mathcal{D}}(\Phi (X),\Phi (Y))$ \textit{is 
injective}. 
\vskip3mm 

\begin{proof}
For $p,i \in {\mathbb Z}$, we endow $\text{T}^pF^iX$ with the filtration 
whose $j$th term is $\text{T}^pF^jX$ for $j > i$ and $\text{T}^pF^iX$ for 
$j\leq i$, and similarly for $\text{T}^pF^iY$. We also endow $\text{T}^pX^i$ 
with the filtration whose $j$th term is $\text{T}^pX^i$ for $j\leq i$ and 0 
for $j > i$, and similarly for $\text{T}^pY^i$. 
We prove (a) and (b) simultaneously, by induction on $N := \text{card}\{
i\in {\mathbb Z}\  \vert \  X^i\neq 0\} + \text{card}\{j\in {\mathbb Z}\  
\vert \  Y^j\neq 0\}$. The case $N\leq 2$ is obvious.

For the {\it induction step}, assume, firstly, that $\text{card}\{j\in 
{\mathbb Z}\  \vert \  Y^j\neq 0\}\geq 2$ and let $n := \text{inf}\{j\in 
{\mathbb Z}\  \vert \  Y^j\neq 0\}$. By applying 
$\text{Hom}_{\mathcal{C}}(X,-)$ to the complex: 
\[
\text{T}^{-1}Y^n\rightarrow F^{n+1}Y\rightarrow Y\rightarrow Y^n
\rightarrow \text{T}F^{n+1}Y 
\tag{*}
\]
and $\text{Hom}_{\mathcal{D}}(\Phi (X),-)$ to the complex $\Phi$((*)), one 
gets a commutative diagram with exact rows and five vertical arrows. If the 
pair $(X,Y)$ verifies the hypothesis of (a) (resp., (b)) then $(X,F^{n+1}Y)$ 
and $(X,Y^n)$ verify the hypothesis of (a) (resp., (b)), and 
$(X,\text{T}F^{n+1}Y)$  verifies the hypothesis of (b) 
(resp., $(X,\text{T}^{-1}Y^n)$ verifies the hypothesis of (a)). 
Using, now, the strong form of the ``Five Lemma'' (see 
\cite{ks1}, Chap. I, Ex 1.8. or \cite{ks2}, (8.3.13)) and taking into account 
the induction hypothesis, one gets the desired conclusion. 

Similarly, if $\text{card}\{i\in {\mathbb Z}\  \vert \  X^i\neq 0\}\geq 2$ 
let $m := \text{inf}\{i\in {\mathbb Z}\  \vert \  X^i\neq 0\}$. One applies 
$\text{Hom}_{\mathcal{C}}(-,Y)$ to the complex: 
\[
\text{T}^{-1}X^m\rightarrow F^{m+1}X\rightarrow X\rightarrow X^m\rightarrow 
\text{T}F^{m+1}X 
\tag{**}
\]
and $\text{Hom}_{\mathcal{D}}(-,\Phi (Y))$ to the complex $\Phi$((**)) and one 
uses again the ``Five Lemma''. 
\end{proof}
\vskip3mm      

Before stating an useful consequence of (B.1.), namely Lemma B.4. below,  
we recall the following well known:

\vskip3mm
{\bf B.2. Lemma.}\quad 
\textit{Let} $X^{\bullet}$ \textit{and} $Y^{\bullet}$ 
\textit{be complexes in an abelian category}  
$\mathcal A$ \textit{and} $n\in {\mathbb Z}$.\\ 
\hspace*{3mm} (a) \textit{If} $X^p=0$ 
(\textit{resp}., $\text{H}^p(X^{\bullet})=0)$ 
\textit{for} $p>n$ \textit{then}: 
\begin{gather*}
\text{Hom}_{\text{K}(\mathcal A)}(X^{\bullet},{\tau}^{\leq n}Y^{\bullet})
\overset{\sim}{\longrightarrow} 
\text{Hom}_{\text{K}(\mathcal A)}(X^{\bullet},Y^{\bullet})\\  
(\textit{resp}.,\  
\text{Hom}_{\text{D}(\mathcal A)}(X^{\bullet},{\tau}^{\leq n}Y^{\bullet})
\overset{\sim}{\longrightarrow}
\text{Hom}_{\text{D}(\mathcal A)}(X^{\bullet},Y^{\bullet})).
\end{gather*}
\hspace*{3mm} (b) \textit{If} $Y^p=0$ 
(\textit{resp}., $\text{H}^p(Y^{\bullet})=0)$ \textit{for} $p<n$ 
\textit{then}:
\begin{gather*}
\text{Hom}_{\text{K}(\mathcal A)}({\tau}^{\geq n}X^{\bullet},Y^{\bullet})
\overset{\sim}{\longrightarrow}
\text{Hom}_{\text{K}(\mathcal A)}(X^{\bullet},Y^{\bullet})\\
(\textit{resp}.,\  
\text{Hom}_{\text{D}(\mathcal A)}({\tau}^{\geq n}X^{\bullet},Y^{\bullet})
\overset{\sim}{\longrightarrow}
\text{Hom}_{\text{D}(\mathcal A)}(X^{\bullet},Y^{\bullet})).
\end{gather*}
\vskip3mm

\begin{proof}
The assertions about $\text{Hom}_{\text{K}(\mathcal A)}$ are easy. 

(a) The inverse of $\text{Hom}_{\text{D}(\mathcal A)}(X^{\bullet},
{\tau}^{\leq n}Y^{\bullet})\rightarrow \text{Hom}_{\text{D}(\mathcal A)}
(X^{\bullet},Y^{\bullet})$ associates to a morphism 
$X^{\bullet}\overset{\text{qis}}{\longleftarrow} X^{\prime \bullet}
\longrightarrow Y^{\bullet}$ in $\text{D}(\mathcal A)$ the morphism 
$X^{\bullet}\overset{\text{qis}}{\longleftarrow} {\tau}^{\leq n}
X^{\prime \bullet}\longrightarrow {\tau}^{\leq n}Y^{\bullet}$. 

(b) The inverse of $\text{Hom}_{\text{D}(\mathcal A)}({\tau}^{\geq n}
X^{\bullet},Y^{\bullet})\rightarrow \text{Hom}_{\text{D}(\mathcal A)}
(X^{\bullet},Y^{\bullet})$ associates to a morphism 
$X^{\bullet}\longrightarrow Y^{\prime \bullet}
\overset{\text{qis}}{\longleftarrow} Y^{\bullet}$ the morphism 
${\tau}^{\geq n}X^{\bullet}\longrightarrow {\tau}^{\geq n}Y^{\prime \bullet} 
\overset{\text{qis}}{\longleftarrow} Y^{\bullet}$.
\end{proof}

\vskip3mm
{\bf B.3. Definition.}\quad 
If $X$, $Y$ are objects of an abelian category $\mathcal A$ and $p\in 
{\mathbb Z}$ then $\text{Ext}^p_{\mathcal A}(X,Y):=
\text{Hom}_{\text{D}({\mathcal A})}(X,\text{T}^pY)$. 
It follows from (B.2.) that $\text{Ext}^p_{\mathcal A}(X,Y)=0$ for $p<0$. 
Moreover, using the arguments from the proof of (B.2.), one sees easily 
that $\text{Hom}_{\mathcal A}(X,Y)\overset{\sim}{\rightarrow} 
\text{Ext}_{\mathcal A}^0(X,Y)$. 
\vskip3mm

The following lemma appears, 
in weaker variants, in several papers like, 
for example, Kapranov \cite{kap} or Canonaco \cite{can}, (A.5.3.). In the 
more precise form (B.4.) below, it was proved in  
\cite{ctr}, (3.3.), under the assumption that the abelian category 
$\mathcal A$ contains sufficiently many injective objects. Here we drop 
this assumption using an argument similar to that used by Canonaco (this 
argument actually appears in the proof of (B.1.)). 

\vskip3mm
{\bf B.4. Lemma.}\quad 
\textit{Let} $\mathcal A$ \textit{be an abelian category}, 
$X^{\bullet}\in 
\text{Ob}\, \text{C}^{-}(\mathcal A)$ \textit{and}  
$Y^{\bullet}\in \text{Ob}\, \text{C}^{+}(\mathcal A)$. 
\textit{Consider the canonical morphism} $\phi :
\text{Hom}_{\text{K}(\mathcal A)}(X^{\bullet},Y^{\bullet})\rightarrow 
\text{Hom}_{\text{D}(\mathcal A)}(X^{\bullet},Y^{\bullet}).$\\
\hspace*{3mm} (a) \textit{If} $\text{Ext}^{p-q}_{\mathcal A}(X^p,Y^q)=0$, 
$\forall p>q$, \textit{then} $\phi $ \textit{is surjective}.\\
\hspace*{3mm} (b) \textit{If} $\text{Ext}^{p-q-1}_{\mathcal A}(X^p,Y^q)=0$, 
$\forall p>q+1$, \textit{then} $\phi $ \textit{is injective}. 
\vskip3mm

\begin{proof}
Let $m:=\text{sup}\{p\in {\mathbb Z}\  \vert \  X^p\neq 0\}$ and 
$n:=\text{inf}\{q\in {\mathbb Z}\  \vert \  Y^q\neq 0\}$. Taking into 
account (B.2.), one may replace $X^{\bullet}$ by ${\tau}^{\geq n}
X^{\bullet}$ and $Y^{\bullet}$ by ${\tau}^{\leq m}Y^{\bullet}$, hence one 
may assume that $X^{\bullet}$ and $Y^{\bullet}$ are {\it bounded} complexes.

In this case, one endows $X^{\bullet}$ with the filtration $F^iX^{\bullet} := 
{\sigma}^{\geq i}X^{\bullet}$ ($\sigma$ = ``stupid truncation''). To the 
semi-split short exact sequence:
\[
0\rightarrow {\sigma}^{\geq i+1}X^{\bullet}\rightarrow 
{\sigma}^{\geq i}X^{\bullet}\rightarrow \text{T}^{-i}X^i\rightarrow 0
\]
one can associate (see, for example, \cite{coa}, (2)(ii)) a distinguished 
triangle in $\text{K}^b(\mathcal{A})$: 
\[
{\sigma}^{\geq i+1}X^{\bullet}\rightarrow {\sigma}^{\geq i}X^{\bullet} 
\rightarrow \text{T}^{-i}X^i\rightarrow \text{T}{\sigma}^{\geq i+1}
X^{\bullet}. 
\]
One also endows $Y^{\bullet}$ with the similar filtration. The conclusion of 
the lemma follows now from (B.1.) applied to the canonical functor 
$\text{K}^b(\mathcal{A})\rightarrow \text{D}^b(\mathcal{A})$. 
The hypotheses of (B.1.) can be easily checked in 
this case because most of the Hom groups involved are zero and 
$\text{Hom}_{\text{K}(\mathcal A)}(X,Y)\overset{\sim}{\rightarrow} 
\text{Hom}_{\text{D}(\mathcal A)}(X,Y)$, $\forall X,Y\in \text{Ob}\, 
{\mathcal A}$ (see (B.3.)).       
\end{proof}


\begin{thebibliography}{30}
\bibitem{abe} T. Abe and M. Yoshinaga, \textit{Splitting criterion for 
              reflexive sheaves}, Proc. Amer. Math. Soc. \textbf{136}, 
              No.6, 1887-1891 (2008). 
\bibitem{bar} D.W. Barnes and L.A. Lambe, \textit{A fixed point 
              approach to homological perturbation theory}, Proc. Amer. 
              Math. Soc. \textbf{112}, No.3, 881-892 (1991).
\bibitem{bgs} A. Beilinson, V. Ginzburg and W. Soergel, \textit{Koszul 
              duality patterns in representation theory}, J. Amer. Math. 
              Soc. \textbf{9}, No.2, 473-527 (1996). 
\bibitem{bgg} I.N. Bernstein, I.M. Gel'fand and S.I. Gel'fand,  
              \textit{Algebraic bundles over} ${\mathbb P}^n$ \textit{and 
              problems of linear algebra}, Funktsional'nyi Analiz i Ego 
              Prilozhenia \textbf{12}, No.3, 66-67 (1978), English 
              translation in: Funct. Anal. Appl. \textbf{12}, 
              212-214 (1978). 
\bibitem{ber} C. Bertone and M. Roggero, \textit{Splitting type, global 
              sections and Chern classes for vector bundles on} 
              ${\mathbb P}^N$, \texttt{arXiv:0804.2985 [math.AG]}.
\bibitem{bro} R. Brown, \textit{The twisted Eilenberg-Zilber theorem},  
              Celebrazioni Archimedee del secolo XX, Simposio di 
              Topologia, 34-37 (1964).
\bibitem{can} A. Canonaco, \textit{The Beilinson complex and canonical 
              rings of irregular surfaces}, Memoirs Amer. Math. Soc. 
              \textbf{862} (2006).
\bibitem{coa} I. Coand\u{a}, \textit{On the Bernstein-Gel'fand-Gel'fand 
              correspondence and a result of Eisenbud, Fl\o ystad, and 
              Schreyer}, J. Math. Kyoto Univ. \textbf{43}, No.2, 
              429-439 (2003).
\bibitem{ctr} I. Coand\u{a} and G. Trautmann, \textit{Horrocks theory 
              and the Bernstein-Gel'fand-Gel'fand correspondence}, 
              Trans. Amer. Math. Soc. \textbf{385}, No.3, 
              1015-1031 (2005).
\bibitem{emc} S. Eilenberg and S. MacLane, \textit{On the group}  
              $\text{H}(\pi ,n)$.I, Ann. Math. \textbf{58}, 
              55-106 (1953). 
\bibitem{eis} D. Eisenbud, \textit{Commutative algebra with a view towards 
              algebraic geometry}, Graduate Texts in Math. \textbf{150}, 
              (Springer-Verlag, 1995). 
\bibitem{efs} D. Eisenbud, G. Fl\o ystad and F.-O. Schreyer, \textit{Sheaf 
              cohomology and free resolutions over exterior algebras},  
              Trans. Amer. Math. Soc. \textbf{355}, No.11, 4397-4426 
              (2003).
\bibitem{gma} S.I. Gel'fand and Yu.I. Manin, \textit{Methods of homological 
              algebra}, (Springer-Verlag, 1996).
\bibitem{gug} V.K.A.M. Gugenheim, \textit{On the chain complex of a fibration}, 
              Illinois J. Math. \textbf{16}, 398-414 
              (1972).
\bibitem{ho1} G. Horrocks, \textit{Vector bundles on the punctured 
              spectrum of a local ring}, Proc. London Math. Soc.  
              \textbf{14}, 689-713 (1964).
\bibitem{ho2} G. Horrocks, \textit{Construction of bundles on} 
              ${\mathbb P}^n$, In: A. Douady and J.-L. Verdier (eds.),  
              \textit{Les equations de Yang-Mills, S\' eminaire E.N.S. 
              (1977-1978)}, Ast\' erisque \textbf{71-72}, 
              (Soc. Math. de France, 1980) 197-203.
\bibitem{kap} M.M. Kapranov, \textit{Derived category of coherent sheaves 
              on Grassmann manifolds}, Funktsional'nyi Analiz i Ego 
              Prilozheniya \textbf{17}, No.2, 78-79 (1983), English 
              translation in: Funct. Anal. Appl. \textbf{17}, 145-146 
              (1983). 
\bibitem{ks1} M. Kashiwara and P. Schapira, \textit{Sheaves on 
              manifolds}, Grundlehren Math. Wiss. \textbf{292}, 
              (Springer-Verlag, 1990). 
\bibitem{ks2} M. Kashiwara and P. Schapira, \textit{Categories and sheaves},  
              Grundlehren Math. Wiss. \textbf{332}, (Springer-Verlag, 2006). 
\bibitem{lst} L. Lambe and J. Stasheff, \textit{Applications of 
              perturbation theory to iterated fibrations}, 
              Manuscripta Math. \textbf{58}, 363-376 (1987). 
\bibitem{ser} J.-P. Serre, \textit{Faisceaux alg\' ebriques coh\' erents},  
              Ann. Math. \textbf{61}, 197-278 (1955).
\bibitem{shi} W. Shih, \textit{Homologie des espaces fibr\' es}, Inst. des 
              Hautes \' Etudes Sci. \textbf{13}, 93-176 (1962).
\bibitem{trm} G. Trautmann, \textit{Moduli of vector bundles on} 
              ${\mathbb P}_n({\mathbb C})$, Math. Ann. \textbf{237}, 
              167-186 (1978).   
\end{thebibliography}
\end{document}